\documentclass[12pt]{amsart}
\usepackage{amsmath, amsfonts, amssymb, amsthm, amscd, graphicx, latexsym, color}
\newcommand{\num}[1]{\# \{ #1 \}}

\newcommand{\norm}[1]{\left\Vert#1\right\Vert}
\newcommand{\abs}[1]{\left\vert#1\right\vert}
\newcommand{\gro}[2]{( #1, #2 )}
\newcommand{\grom}[3]{ ( #1, #2  )_#3}

\newcommand{\HH}{\mathcal{H}}

\newcommand{\lab}[1]{ \label {#1}  }

\newcommand{\bW}{\overline{W}}
\newcommand{\bX}{\overline{X}}

\newcommand{\bx}{\overline{x}}
\newcommand{\bt}{\overline{t}}
\newcommand{\bm}{\overline{m}}
\newcommand{\bo}{\overline{o}}
\newcommand{\ov}[1]{\overline{#1}}

\newcommand{\hyp}{$\delta$--hyperbolic}

\newtheorem{thm}{Theorem}[section]

\newtheorem{prop}[thm]{Proposition}
\newtheorem{coro}[thm]{Corollary}
\newtheorem{lemm}[thm]{Lemma}
\newtheorem{rem}[thm]{Remark}

\newtheorem{defn}[thm]{Definition}

\newcommand{\ga}{\gamma}

\newcommand{\eps}{\epsilon }

\newcommand{\cay}{$C a  y (G)$}

\begin{document}
\title[Actions of maximal growth]{Actions of maximal growth of hyperbolic groups}

\author{Vladimir Chaynikov}
\maketitle
\begin{abstract} We prove that  every non-elementary hyperbolic group $G$ acts with maximal growth on some set $X$ such that every orbit of any element $g \in G$ is finite. As a side-product of our approach we prove that for a  non-elementary hyperbolic group  $G$ and a quasiconvex subgroup of infinite index $\HH \leq G$ there exists $g \in G$ such that $\langle \HH,g \rangle$ is quasiconvex of infinite index and is isomorphic to $\HH*\langle g \rangle$ if and only if $\HH \cap E(G)= \{ e \} $, where $E(G)$ is the maximal finite normal subgroup of $G$. 
\end{abstract}

\subjclass{AMS2001:}  {20F67}

\section{Introduction}
The notion of growth of algebraic structures has been extensively studied. In the case of groups, there are three main classes of growth: polynomial, subexponential, exponential.  In \cite{BO} the authors discuss the notion of growth of actions   of a group (monoid, ring) on a set (vector space). Let us denote the growth function of a transitive action of a group $G$ generated by $S$ on a set $X$ with respect to some base point $o \in X$ by 
$g_{o,S}(n)=\num{o'=og \vert \abs{g} \leq n   }$ (see section \ref{Actions of Maximal Growth}). 

A distinguished class of actions defined and studied in \cite{BO} is that of actions of maximal growth. We observe that in case of a non-amenable group $G$ the growth of action of $G$ on $X$ is \emph{maximal} if there exists $c_1>0$ such that 
$$c_1 f(n) \leq g_{o,S}(n)$$
for every natural $n$, where $f(n)$ is the growth of the group $G$ itself (see remark \ref{morethenexponent}).

In \cite{BO} the authors construct some examples of actions by the free group of maximal growth and satisfying additional properties, see for example corollary \ref{coro6}. The main result of our paper is the following  broadening of the aforementioned corollary: 

\begin{thm} \lab{maxgrowth_main} Let $G$ be a non-elementary hyperbolic group. Then there exist a set $X$ and an action of $G$ on $X$ such that the growth of this action is maximal and each orbit of action by every element $g \in G$ is finite.
\end{thm}
One can observe that the  above result follows  from theorem  \ref{maxgrowthburnside} in this paper.

The main theorem stems from the technical result (theorem) \ref{quasiconvex and isomorphic} which also allows us to generalize and strengthen the result of Arzhantseva \cite{Arzh} conjectured by M. Gromov \cite{Gro}.

The following  theorem and corollary generalize theorem 1 in \cite{Arzh} by removing the requirement on the hyperbolic group to be torsion-free  and formulating  necessary and sufficient conditions.  We recall the notation $E(g)$ -- the maximal elementary subgroup of hyperbolic group $G$ containing $g$ (it exists whenever $g$ is of infinite order, see section \ref{hypcase}). Recall also that there exists a unique finite maximal normal subgroup $E(G)$ in every non-elementary hyperbolic group $G$. We will call $E(G)$ \emph{the finite radical}\footnote{the term proposed by A. Olshanskiy.} of $G$.
\begin{thm} \lab{generalizedA3}  Let $G$ be a non-elementary hyperbolic group and $\HH$ be a quasiconvex infinite index subgroup of $G$. 

(1) Consider an element $x$ in $G$ of infinite order.  
There exists  a natural number  $ t>0$ such that the subgroup  $\langle \HH, x^t \rangle$ (i.e. generated by $\HH$ and $x^t$) is  isomorphic to $\HH*\langle  x^t \rangle$ if and only if $E(x) \cap \HH=\{e \}$. 

(2) \footnote{When preparing this paper for  publication the author learned that  a version of this statement have been presented by F. Dudkin and K. Sviridov on a Group Theory seminar in IM SORAN (Novosibirsk) in November, 2011.} An element $x$, satisfying part (1), exists if and only if $\HH \cap E(G)=\{e \}$.

(3) If $\HH \cap E(G)=\{e \}$ then for $x$ and $t$ described in part (1) the subgroup $\langle \HH, x^t \rangle$ is quasiconvex of infinite index and the intersection $E(G) \cap \langle \HH, x^t \rangle$ is  trivial.
\end{thm}

Part (1) of theorem \ref{generalizedA3} follows also from a more general statement in \cite{M-P}(Corollary 1.12) and a particular case when $E(x)=E^+(x)$ appears in theorem 5 \cite{Min}.
We also formulate the following (somewhat more general) result concerning  arbitrary quasiconvex subgroups of infinite index.  
 
\begin{coro} \lab{generalizedA2} Let $G$ be a non-elementary hyperbolic group and $\HH$ be a quasiconvex subgroup of infinite index in $G$. Then there exists $g \in G$ of infinite order such that  $\langle \HH \cdot E(G), g \rangle \cong \HH \cdot E(G)  * _{E(G)}\langle g,E(G) \rangle. $ Moreover $\langle \HH \cdot E(G), g \rangle$ is a  quasiconvex subgroup of infinite index.  
 \end{coro}

\section{Hyperbolic spaces and hyperbolic groups} \lab{hypcase}
 {\bf Hyperbolic Spaces}. We recall some definitions and properties from the founding article of Gromov \cite{Gro} (see also \cite{Ghys}).
 Let $(X, \abs{\ })$ be a metric space. We sometimes denote the distance $\abs{x-y}$  between $x,y \in X$ by $d(x,y)$. We assume  that $X$ is geodesic, i.e. every two points can be connected by a geodesic line. We refer to a geodesic between some point $x,y $ of $X$ as $[x,y]$. For convenience we denote by $\abs{x}$ the  distance $\abs{x-y_0}$ to some fixed point $y_0$ (usually the identity element of the group). 
 
 For a path $\ga$ in $X$ we denote the initial (terminal) vertex of $\ga$ by $\ga_-$ ($\ga_+$), denote by $\norm \ga$  the length of path $\ga$ and by $\abs \ga$  the distance $\abs{\ga_+ - \ga_-}$. Recall that if $0<\lambda \leq 1$ and $c\geq 0$ then a path $\ga$ in $X$ is called \emph{$(\lambda,c)$--quasigeodesic} if for every subpath $\ga_1$  of $ \ga$ the following inequality is satisfied:
 
 $ \norm {\ga_1} \leq \frac{1} {\lambda} \abs{\ga_1} +c$.

 We call the path $\ga$  \emph{geodesic up to $c$}, if it is $(1,c)$-quasigeodesic.

Define a scalar (Gromov) product of $x,y$ with respect to $z$ by formula
 $$ \grom{ x}{y} {z}=\frac{1}{2} (\abs{x-z}+\abs{y-z}-\abs{x-y}).$$

An (equivalent) implicit definition of the Gromov product illustrates its geometric significance:
\begin{equation} \lab{Gromprod} \grom{x}{y}{z}+\grom{x}{z}{y}=\abs{z-y};
\end{equation}
\begin{equation*}\grom{x}{y}{z}+\grom{y}{z}{x}=\abs{z-x};
\end{equation*}
\begin{equation*} \grom{y}{z}{x}+\grom{z}{x}{y}=\abs{x-y}.
\end{equation*}

A space $X$ is called \hyp \  if there exists a non-negative integer $\delta$ such that the following  inequality holds:
$$ {\bf (H1)} \  \  \forall x,y,z,t \in X, \ \  \grom{ x} {y} {z} \geq min\{ \grom{ x} {t} {z}, \grom{ y}{t}{z} \} - \delta.  $$  
 
The condition (H1) implies (and in fact is equivalent up to constant) the following:

\emph{ {\bf(H2)} For every triple of points $x,y,z$ in $X$ every  geodesic $[x,y]$ is within the (closed) $4\delta$-neighborhood of the union $[x,z] \cup [y,z]$. } 
 
\emph{ {\bf(H3)} For every four points $x,y,z, t$ in $X$ we have $\abs{x-y}+\abs{z-t} \leq max \{ \abs{x-z}+\abs{y-t}, \abs{x-t}+\abs{y-z}  \}+2\delta.$  }

 



 



\medskip {\bf Hyperbolic Groups.} Let $G$ be a finitely presented group with presentation $gp ( S \vert \mathcal{D} )$. We assume that no generator in $S$ is equal to $e$ in $G$. We consider $G$ as a metric space with respect to the distance function 
$\abs{g-h}=\abs{g^{-1}h}$ for every $g$ and $h$. We denote by $\abs g$  the length of a minimal (geodesic) word  with respect to the generators $S$ equal to $g$. The notation  $\gro{g}{h} $ is the Gromov product $\grom{g} {h}{e}$ with respect to the identity vertex $e$.

We denote the (right) Cayley graph of the group by  $Cay(G)$. The graph $Cay(G)$ has a set of vertices $G$,  and a pair of vertices $g_1,g_2$ is connected by an edge of length $1$ labeled by $s$ if and only if $g^{-1} _1 g_2=s$ in $G$ for some $s \in S^{\pm 1}$. It is clear that \cay \ may be considered  as a geodesic space: one identifies every edge of \cay \ with interval $[0,1]$ and chooses the maximal metric $d$ which agrees with metric on every edge.
Define a label function on paths in $Cay(G)$. From now on, by a path in $Cay(G)$ we mean a path $p=p_1 ... p_n$, where  $p_i$ is an edge in \cay \ between some group elements $g_i$, $g_{i+1}$ for every $1\leq i \leq n$. A label $lab (p)$ function is defined on any path $p$    by equality $lab(p)=lab(p_1) ... lab(p_n)$, i.e. $lab(p)$ is a word in alphabet $S^{\pm 1}$.

Hence a  unique word $lab(p)$ is assigned to a path $p$ in $Cay(G)$. On the other hand for every word $w$ in alphabet $S^{\pm 1}$ there exists a unique path $p$ in $Cay(G)$ starting from the identity vertex with label $w$. Hence there is a one-to-one correspondence between paths with initial vertex $e$ (the identity vertex in $G$) and words in alphabet $S^{\pm 1}$, so we will not distinguish between a word in the alphabet $S^{\pm 1}$ and it's image in $Cay(G)$, i.e. a path starting from the identity vertex. Thus, when considering some words $X,Y,Z$  in the alphabet  $S^{\pm 1 }$, we  can talk about the path $\gamma=XYZ$ in the Cayley graph of $G$ originating in the identity vertex $e$. To distinguish a path $Y$ with initial vertex $e$ from the subpath of $\ga$ with label $Y$ we denote the latter as $_\ga Y$. 
We will talk about values $\abs X, \norm X$ for a word $X$ in alphabet  $S^{\pm1}$ meaning these values on the corresponding paths in \cay. Given elements $x_1, ... , x_k$ in $G$ we may write $lab(p)= x_1^{t_1}...x_k^{t_k}$ for some path $p$ in \cay, $t_i \in \mathbb{Z}$ if for some geodesic words $X_1,...,X_k$ representing elements $x_1,...,x_k$ we have  $lab(p) = X_1^{t_1}...X_k^{t_k}$.

  For a point $x$ in a metric space $X$ and $r \geq 0$ we denote by $B_r (x)$ a metric ball of radius $r$ around $x$. For a set $D \subset X$ we denote by $B_r(D)$ a (closed) $r$-neighborhood of $D$ in $X$ (i.e. $B_r(D)=\cup_{x \in D} B_r(x)$). We denote the ball $B_R(e)$ in the  Cayley graph \cay \ by $B_R$. Given a set $D \in Cay(G)$ we denote by $\num {D}$ a number of vertices in $D$.

A group $G$ is called \emph{$\delta$-hyperbolic} for some $\delta \geq 0$, if it's Cayley graph is $\delta$-hyperbolic. It is well known that hyperbolicity does not depend on choice of a finite presentation of the group $G$ (while $\delta$ does depend on presentation).

\begin{lemm} \lab{lemma 1.9} (\cite{Gro}, \cite{Ghys} p. 87) There exists a constant $H=H(\delta,\lambda,c)$ such that for any $(\lambda,c)$--quasigeodesic path $p$ in a \hyp \ space and any geodesic path $q$ with conditions $q_-=p_-$ and $q_+=p_+$, the  paths $p$ and $q$ are within (closed) $H$-neighborhoods of each other.
\end{lemm}

We recall that a (sub)group is called elementary if it contains a cyclic group of finite index. For any element $g \in G$ of infinite order in a hyperbolic group, there exists a unique maximal elementary subgroup $E(g)$  containing $g$ (see \cite{Gro}, \cite{Olsh93} lemma 1.16). It is well known that for a hyperbolic group $G$
$$E(g)=\{x \in G \vert \exists n\neq 0 \text{ such that } xg^nx^{-1}=g^{\pm n} \text{ in $G$} \}, $$
and if $a$ is an element in  $ E(g)$ of infinite order then $E(g)=E(a)$. We recall also that if $G$ is a non-elementary hyperbolic then the  subgroup $E(G)=\cap_g \{E(g)  \vert g \in G,\text { order of $g$ is infinite} \}$ is a unique maximal finite normal subgroup (\cite{Olsh93}, prop.1). As agreed in the introduction, we will call $E(G)$ the \emph{finite radical} of a non-elementary group $G$. 

\begin{defn} A subset $A$ is called $K$-quasiconvex in the metric space $X$  if for any pair of points $a,b \in A$ every geodesic connecting $a$ and $b$ (in $X$) is within (closed) $K$-neighborhood of $A$.   A subgroup $\HH$ of a hyperbolic group $G$ is $K$-quasiconvex if it forms a  $K$-quasiconvex subset in the graph \cay.
\end{defn}
It is said that $\HH$ is quasiconvex if it is $K$-quasiconvex for some $K \geq 0$.
Note also that the left multiplication $g \rightarrow a g$ induces an isometry of $G$ and hence, for a $K$-quasiconvex subgroup $\HH$, the right coset $a\HH$ is $K$-quasiconvex  for any $a$ in $G$. 

\begin{lemm} \lab{width} (\cite{GMRS}, lemma 1.2) Let $H$ be a $K$-quasiconvex subgroup of a \hyp \ group $G$. If the shortest representative of a double coset $HgH$ has length greater than $2K+2\delta$, then the intersection $H \cap g^{-1}Hg$ consists of elements shorter than $2K+8\delta+2$ and, hence, is finite. 
\end{lemm}

\begin{prop} \lab{doublecoset} (\cite{Arzh}, prop.1) Let $G$ be a word hyperbolic group and $H$ a quasiconvex subgroup of $G$ of infinite index. Then the number of double cosets of $G$ modulo $H$ is infinite.  
\end{prop}

We quote the following:
\begin{thm} \lab{Mack's} (\cite{Mack}, theorem 6.4) Let $G$ be a  hyperbolic group and $H$ be a quasiconvex subgroup of infinite index in $G$. Then there exist $C>0$ and a set-theoretic section $s:G/H \rightarrow G$ such that: 

(1) the section $s$ maps each coset $gH$ to an element $g' \in gH$ with $\abs{g'} $  minimal among all representatives in $gH$; 

(2) the group $G$ is within $C$-neighborhood of $s(G/H)$.
\end{thm}

The following lemma summarize some properties of elementary subgroups of hyperbolic groups (\cite{Gro}; \cite{Ghys} p.150, p.154; \cite{CDP} Pr. 4.2, Ch.10; \cite{Olsh93} lemma 2.2).

\begin{lemm} \lab{lemma 1.11} Let $G$ be a hyperbolic group. 

(i) For any word $W$ of infinite order in the hyperbolic group $G$ there exist constants $0< \lambda \leq 1$ and $c\geq 0 $ such that any path with label $W^m$ in \cay \ is  $(\lambda,c)$--quasigeodesic  for any $m$.

(ii)  Let $E$ be an infinite elementary subgroup in $G$. Then there exists a constant $K=K(E) \geq 0$ such that the subgroup $E$
is $K$-quasiconvex.

(iii) If $W$ is a geodesic word and $p$ is a path with label $W^n$ then there exists $K$ (independent of $n$) such that the path $p$ and the geodesic $[p_-,p_+]$ are within $K$--neighborhoods of each other. 
 
(vi) Let $g,h$ be elements of infinite order such that $E(g)\neq E(h)$. Then the Gromov products $\gro{g^m}{h^n}, \gro{g^u}{g^v}, \gro{h^u}{h^v}$ are bounded by some constant $C$ depending on $g,h$ only provided $uv<0$.
 \end{lemm}
Following \cite{Olsh93}, we call a pair of elements $x, y$ of infinite order in $G$ non-commensurable
if $x^k$ is not conjugate to $y^s$ for any non-zero integers $k, s$.

\begin{lemm} (\cite{Olsh93}, lemmas 3.4, 3.8)\lab{3.4,3.8} There exist infinitely many pairwise non-commensurable elements $g_1,g_2, ...$ in a non-elementary hyperbolic group $G$ such that $E(g_i)=\langle g_i \rangle\times E(G)$ for every $i$.
\end{lemm}

 Let $W$ be a word, and let us fix some factorization $W \equiv W_1^{i_1}W_2^{i_2} ... W_k^{i_k}$ for some  words $W_1, ..., W_k$. Consider a path $q$ with label $W$ in \cay.
 
 Consider all vertices $o_i$ which are the terminal vertices of initial subpaths $p_i$ of $q$ such that $lab(p_i)= W_1^{i_1}...W_{m-1}^{i_{m-1}}W_m^s$, where $m \leq k $ and $s=0,...,i_m$. Following \cite{Olsh93}, we call vertices $\{ o_i \}$ \emph{phase vertices} of $q$ relative to factorization  $W_1^{i_1}W_2^{i_2} ... W_k^{i_k}$ of the $lab(q)$. We enumerate distinct phase vertices along the path $q$ starting from $o_0=q_-$; the total number of such vertices is $(\abs{i_1}+... +\abs{i_k}+1)$.  
 
Assume we have a pair of paths $q,\bar{q}$ in \cay \  with phase vertices $o_i$ and $\bar{o}_j$ where
 $i=1, ..., l$, $j=1, ... , m $ for some positive integers $l,m$. We call a shortest path between a phase vertex $o_i$ and some phase vertex $\bar{o}_j$ of $\bar{q}$ a \emph{phase path} with initial vertex $o_i$. We may also talk about phase vertices of subpaths $p$ of $q$ meaning these vertices $o_i$ which belong to $p$.

\begin{defn} \cite{Olsh93}  Let the words $W_1, ..., W_l$   represent some elements of infinite order in $G$.  Fix  some $A \geq 0$ and  an integer $m$ to define a set $S_m=S(W_1, ...,W_l, A, m)$  of words 
$$W=X_0W_1^{m_1}X_1 W_2^{m_2} ... W_l^{m_l}X_l \text{ where } \abs{m_2}, ..., \abs{m_{l-1}}\geq m,$$ 
 such that $\norm {X_i} \leq A$ for $i=0,..., l$ and $X_i^{-1} W_iX_i \notin E(W_{i+1})$ in $G$ for $i=1, ... ,l-1$. If $l=1$ we assume that $\abs{m_1} \geq m$.
\end{defn}
  
 \begin{lemm} (\cite{Olsh93}, lemma 2.4)\lab{O2.4} There exist $\lambda>0$, $c \geq 0$ and $m>0$ (depending on $K,W_1,W_2,..., W_l$) such that  any word $W \in S_m$ is $(\lambda,c)$--quasigeodesic. If $W_i \equiv W_j$ for all $i,j$ then the constant $\lambda$ does not depend on $A,l$. 
\end{lemm} 

  Consider a closed path $p_1q_1p_2q_2$ in \cay. Let $q_1=x_1t_1x_2t_2 ... x_lt_l$ where $lab(x_i)= X_i$ and $lab( t_i)= W_i^{m_i}$ for some 
  $W=X_0W_1^{m_1}X_1 W_2^{m_2} ... W_l^{m_l}X_l \in S_m$. Similarly, we let $q_2^{-1}=\bx_1 \bt_1 ... \bx_l \bt_l$ where $lab(\bx_i)= \bX_i$ and $lab( \bt_i)= \bW_i^{\bm_i}$ for some   $\bW=\bX_0 \bW_1^{\bm_1}\bX_1 \bW_2^{\bm_2} ... \bW_l^{\bm_l}\bX_l \in \ov S _m$. Define phase vertices $o_i$ and $\bo_j$ on $q_1$ and $q_2$ relative to factorizations $X_0W_1^{m_1}X_1 W_2^{m_2} ... W_l^{m_l}X_l$ and $\bX_0 \bW_1^{\bm_1}\bX_1 \bW_2^{\bm_2} ... \bW_l^{\bm_l}\bX_l$. As in \cite{Olsh93}, We say that  paths $t_i$ and $\bt_j$ are compatible if there exists a phase path $v_i$ with $lab(v_i)=V_i$ between a phase vertex of $t_i$ and $\bt_j$ such that there exist natural numbers $a,b$ satisfying  $(V_i\bW_jV_i^{-1})^a=W_i^b$. 
  
\begin{lemm} (\cite{Olsh93}, lemma 2.5)\lab{lemma 2.5} 
Provided the conditions for $q_1$ and $q_2$ hold, and $\abs{p_1}, \abs{p_2}<C$ for some $C$,  there exists an integer $m$ and an integer $k$, where $\abs k \leq 1$ such that $t_i$ and $\bar{t}_{i+k}$ are compatible for any $i=2,...,l-1$ provided that  $\abs{m_2}, ... ,\abs{m_{l-1}},\abs{\bm_2}, ... ,\abs{\bm_{l-1}}\geq m$ and for  $i=1$ (resp. $ i=l$) if  $\abs{m_1}\geq m,$ (resp. $ \abs{m_l}\geq m$). Moreover $t_i$ is not compatible with $\bar t_j$ if $j \neq i+k$.
\end{lemm}   


\section{Actions of Maximal Growth} \lab{Actions of Maximal Growth}
Let $G$ be a group generated by a finite set $S$ and suppose that $G$ acts on a set $X$ from the right:
$$xe=x, \ \ (xg_1)g_2=x(g_1g_2) \text{ in  $G$ for all } x \in X; \ g_1,g_2 \in G.$$
We assume that the action is transitive (i.e. $X=oG$, where $o$ is some element from $X$). Consider the set $B_{n}(o)$ of elements $og \in X$ such that $g \in G$ and $\abs{g} \leq n$. Then \emph{the growth function of the right action } of $G$ on $X$ is $f_{o,S} (n)=\num{B_{n}(o)}$. 
Let $o'=og_0 \in G$ and denote $\abs{g_0}$ by $C$. It is clear that $B_{n}(o') \subset B_{n+C}(o)$ and hence 
$$f_{o,S}(n+C) \geq f_{o',S}(n).$$

Consider a set $\mathcal{F}$ of functions from $\mathbb N _0$ to $\mathbb N _0$.  A pair $f,g \in \mathcal F$ is said (see \cite{BO}, \S 1.4) to satisfy the relation $f \prec g$ if there exist a non-negative integer $C$ such that $f(n) \leq g(n+C)$ for every $n \in \mathbb N _0$. Clearly the relation $ \prec$ is transitive and reflexive.  Functions $f,g \in \mathcal F$ are said to be \emph{equivalent} (\cite{BO}, \S 1.4)  if $f \prec g$ and $g \prec f$. According to the discussion above, growth functions of transitive action of $G$ on $X$ with respect to different base points $o,o'$ are equivalent. 

  If the group $G$ acts from the right on $X=G$,  we get the usual growth function and denote it by $f(n)$; clearly the growth of any action of $G$ is bounded by the usual growth function of $G$: $f_{o,S} (n) \leq f(n)$ for any $o \in X$. 

If $H$ is a stabilizer of $o$, then every element $x \in X$ is in one-to-one correspondence with a coset $Hg$ in $G$ such that $x=og$ and the right actions of $G$ on $X$ and on  $H \backslash G$ are isomorphic.

\begin{defn} (\cite{BO}, \S 2) Let $f(n)$ be a growth function of $G$ relative to a finite generating set $S$ and consider a transitive action of $G$ on a set $X$. Then the growth of the action  is called maximal if the function $f_{o,S}(n) $ is equivalent to $f(n)$.  
\end{defn}

In this paper we discuss the growth of actions of hyperbolic groups which are known to be non-amenable (see remark \ref{hypisnonamenable}).
We recall that a group $G$ is called amenable if there exists a finitely additive left invariant probability measure on $G$ (see \cite{Gre}).  

\begin{rem} \lab{hypisnonamenable}
Every non-elementary hyperbolic group is non-amenable.
\end{rem} 
{\bf Proof} If $G$ is a non-elementary hyperbolic group, then it contains a free non-cyclic subgroup $F_2$ (\cite{Gro}, \cite{Ghys}, p. 157). But a free group of rank greater then one is non-amenable (see \cite{Gre} 1.2.8). On the other hand, a subgroup of an amenable group is amenable (\cite{Gre}, theorem 1.2.5). Hence $G$ cannot be amenable.$\Box$ 

 The famous F\"olner amenability criterion (\cite{Gre}) yields the following: 
\begin{coro} \lab{amenable} For every non-elementary hyperbolic group $G$ there exists $\eps>0$ (depending on $G$ only) such that $\num{B_{R+1}}  \geq (1+\eps)  \num{B_{R}}$ for any $R$.
\end{coro} 

\begin{rem} \lab{morethenexponent} Let $G$ be a non-amenable group with growth function $f(x)$ relative to a finite generating set $S$. Assume $G$ acts on $X$ with respect to some base point $o \in X$; denote the growth function of this action  by $g_{o,S}(x)$. Then  there exists $c_1>0$ such that the inequality $g_{o,S}(n) \geq c_1 f(n)$ holds for all natural $n $ if and only if the action has maximal growth.  
\end{rem}
{\bf Proof} We first show the "only if" part. By corollary \ref{amenable} there exists $\eps>0$ such that the recursive formula $f(n+1) \geq(1+\eps) f(n)$ holds for every $n$.  We choose a natural $C$ satisfying $c_1(1+\eps)^C \geq 1.$ Then, applying the recursive formula $C$ times we get:
$$g_{o,S}(n+C) \geq c_1f(n+C) \geq c_1(1+\eps)^Cf(n) \geq f(n) .$$

Now assume that the action has maximal growth, i.e. $g_{o,S}(n+C) \geq f(n)$ for some natural number  $C\geq 0$ and every natural $n$. It is clear from definition of $g_{o,S}$ that  $ g_{o,S}(n+C) \leq (2\num S)^C\times g_{o,S}(n)$ and hence for $c_1=  (2\num S)^{-C}$ the inequality $g_{o,S}(n) \geq c_1 f(n)$ holds.$\Box$

We recall the notion of exponential growth rate of a group $G$ with respect to the set of generators $S$:
$ \lambda (G,S)=lim_{n \rightarrow \infty} \sqrt[n]{f(n)},$ where $f(n)$ is a growth function of $G$. 

Let $S$ be a finite generating set in $G$ and let $N$ be an infinite normal subgroup of $G$. We denote the image of $S$ under the canonical homomorphism  $G \rightarrow G/N$ by $\overline S$. The following theorem is often summarized by saying that the hyperbolic groups are "growth tight":

\begin{thm} \lab{AL1}\cite{AL} Let $G$ be a non-elementary hyperbolic group and $S$ any finite set of generators for $G$. Then for any infinite normal subgroup $N$ of $G$ we have $\lambda(G,S)>\lambda(G/N,\overline {S})$.
\end{thm}

The next corollary restates the above theorem in terms of maximal growth.



\begin{coro} 
Assume $G$ is a hyperbolic group acting on some set $X$ from the right  with maximal growth. Then the kernel of this action  is a finite normal subgroup.
\end{coro}
{\bf Proof }
Let $N$ be the kernel of the action on $X$. For any point $o \in X$ we have that $oNg=og$ for all $g \in G$ and hence 
the growth function of the action $g_{o,S}$   satisfies:  
\begin{equation} \lab{finite1}
g_{o,S}(n) \leq f_{G/N}(n) \text{ for every } n \in \mathbb N,
\end{equation}
where $\ov f(n)$ is the growth function of $G/N$ with respect to images $\overline {S}$ of generators $S$ of $G$. If $f(n)$ is the growth function of $G$ with respect to $S$ and  the  growth of the action is maximal, then there exists $c_1>0$ such that $g_{o,S}(n) \geq c_1 f(n)$ for for every $ n \in \mathbb N$. Hence, using (\ref{finite1}) and the last inequality, we have:
$$\lambda(G/N,\overline {S})=lim_{n \rightarrow \infty} \sqrt[n]{\ov f (n)} \geq  lim_{n \rightarrow \infty} \sqrt[n]{g_{o,S}(n)} \geq lim_{n \rightarrow \infty} \sqrt[n]{c_1f(n)}=\lambda(G,S),$$
which by Theorem \ref{AL1} can only hold when $N$ is finite. $\Box$


In \cite{BO} the authors provide examples of maximal growth actions of free groups satisfying some additional properties. 

Recall that in \cite{Sta}  a subgroup $H$ of a group $G$ is said to satisfy the Burnside condition if for any $a \in G$ there exists a natural number $n \neq 0$ such that $a^n$ is in $H$.  One of the main results of the aforementioned paper is the following: 

\begin{coro} (\cite{BO}, corollary 6) \lab{coro6} 
Any finitely generated subgroup $H$ of infinite index in the free group $F$ of rank greater then $1$ is contained as a free factor in a free subgroup $K$ satisfying the Burnside condition. One can choose $K$ with  maximal growth of action of $F$ on $K\backslash F$. It follows that there exists a transitive action of $F$, with maximal growth and with finite orbits for each element $g\in F$, which factors through the action of $F$ on $H\backslash F$. 
\end{coro}

The following theorem generalizes the corollary 6 of \cite{BO} from free groups to non-elementary hyperbolic ones:
\begin{thm} \lab{maxgrowthburnside} Let $G$ be a non-elementary hyperbolic group with growth function $f(n)$. Then for any $0<q<1$ there exists a free subgroup $H$ in $G$ satisfying the Burnside condition and such that the growth $f_{H \backslash G}(n)$ of right action of $G$ on $H \backslash G$ satisfies 
$f_{H \backslash G}(n) \geq q f(n)$. In particular, the growth of such action is maximal.
\end{thm}
Consequently, for every non-elementary hyperbolic group $G$ there exists  a transitive action of $G$ with maximal growth such that the orbit of action of any element $g \in G$ is finite.

Throughout this paper we will mainly discuss properties of left cosets. The connection between the right and left cosets is established by the following observation: 
\begin{rem} \lab{symmetry}
 The right coset $Hg$ intersects the ball $B_R$ in $G$ if and only if the left coset $g^{-1}H$ intersects $B_R$.$\Box$
 \end{rem}  
The abundance  of actions of maximal growth is evident from the  following:
\begin{coro} \lab{coroM}  Let $G$ be a  hyperbolic group and $H$ be a quasiconvex subgroup of infinite index in $G$. Then the natural right action of $G$ on $H\backslash G$ has maximal growth. 
\end{coro}
{\bf Proof}  We  first consider  left cosets $G/H$. By  theorem \ref{Mack's} there exists $C >0$ and the section $s$ such that the group $G$ is in $B_C(s(G/H))$. Hence for every $g \in B_R$ there exists $\ov g \in s(G/H)$ such that $\abs{g-\ov g } \leq C$. By definition of $s$, $\abs{ \ov g } \leq \abs g$ and thus $\ov g \in B_R$. We get that $B_R \subset \cup_{g \in B_R \cap s(G/H)} B_C(g)$, which implies:
\begin{equation} \lab{coroM1}
f(R)=\num{B_R} \leq \num{B_C} \times \num {s(G/H) \cap B_R}.
\end{equation}
 If $g_1,g_2 \in s(G/H) \cap B_R$ then (because the map $s$ is a section) $g_1H \neq g_2H$. We get that $\num{s(G/H) \cap B_R} \leq \num {gH \vert gH \cap B_R \neq \emptyset}$ and  the remark \ref{symmetry} provides $\num {gH \vert gH \cap B_R \neq \emptyset}=\num {Hg^{-1} \vert Hg^{-1} \cap B_R \neq \emptyset}$, thus
\begin{equation} \lab{coroM2}
\num{s(G/H) \cap B_R} \leq \num {Hg^{-1} \vert Hg^{-1} \cap B_R \neq \emptyset}. 
\end{equation}
 Evidently the sets $\{ Hg^{-1} \vert Hg^{-1} \cap B_R \neq \emptyset \}$ 
  and $\{ Hg \vert \exists g_1\  : \ Hg_1 = Hg \ \& \  \abs{g_1} \leq R \}$ 
  contain the same cosets, and by definition of the growth function  $f_{H,G/H} (R)  $ of natural right action of $G$ on  $G/H$: $\num {Hg \vert \exists g_1  :  Hg_1 = Hg  \&   \abs{g_1} \leq R}=f_{H,G/H} (R).$  Using inequalities (\ref{coroM1}), (\ref{coroM2}) and the last equality we get:
$$f(R) \leq  \num{B_C} \times \num {s(G/H) \cap B_R} \leq \num{B_C} \times \num{Hg \vert \exists g_1  :  Hg_1 = Hg \  \&   \abs{g_1} \leq R} $$
$$\leq \num{B_C} \times f_{H,G/H} (R).$$
By remark \ref{morethenexponent} the action has maximal growth.$\Box$

 
\section{Proof of theorem \ref{maxgrowthburnside} }

Throughout this paragraph we assume that the group $G$ is non-elementary hyperbolic.  
The following lemma summarizes some geometric properties that we will need later.  
\begin{lemm} \lab{lemma1*} Let $a,b,c,d$ be points in a \hyp \ space $X$.

(i)Assume that $\grom{a}{c}{b}, \grom{b}{d}{c} \leq M$. If we take $Q$ such  that $\grom{a}{d}{b} \leq Q$ then
 $\abs{a-d} \geq \abs{a-b} +\abs{b-c}+ \abs{c-d} -2M-2Q$. Moreover, if $\abs{b-c}>2M+\delta$, then we can choose $Q=M+\delta$.
 
 Assume that the point $d$ is on the segment $[a,b]$ and 
 
(ii) $d\in B_{M_1}([a,c])$ for some $M_1 \geq 0$. Then  
\begin{equation} \lab{lemma1*1}
\abs{d-b}\geq \grom{a}{c}{b}-\delta -M_1 \text{ and}
\end{equation} 
\begin{equation} \lab{lemma1*2}
\abs{a-c}\geq \abs{a-b}+ \abs{b-c} -2\abs{d-b}-2\delta-2M_1.
\end{equation}

(iii)  that  $\grom{b}{c}{a}-5\delta>\abs{a-d}$. Then $d\in B_{4\delta}([a,c]).$

(iv) the vertex $d$ is at least $D>0$ away from each point $a,b$. Then $$\abs{d-c} \leq max \{ \abs{a-c},\abs{b-c} \}+2\delta -D.$$  
\end{lemm}   
{\bf Proof} 
(i) By definition of Gromov product and conditions of part (i) we get that 
$$\abs{a-d}=\abs{a-b}+\abs{b-d}-2\grom{a}{d}{b} = \abs{a-b}+(\abs{b-c} +\ \abs{c-d}-2\grom{b}{d}{c}-2\grom{a}{d}{b}) \geq$$
$$\geq \abs{a-b}+ \abs{b-c}+\abs{c-d}-2M-2Q.$$

It remains to show the second claim in part (i). If $\abs{b-c}>2M+\delta$, then by (\ref{Gromprod}):
\begin{equation} \label {1*1}
\grom{c}{d}{b}=\abs{b-c}-\grom{b}{d}{c}>2M+\delta-M.
\end{equation}

 By inequality (\ref{1*1}) and definition (H1) of \hyp \ space we have 
  $$\grom{c}{d}{b} >M+\delta \geq \grom{a}{c}{b}+\delta \geq  min \{\grom{a}{d}{b}, \grom {c}{d}{b}\},$$ 
which implies that $\grom{a}{d}{b}\leq M+\delta$. 

(ii) Let $d'$ be a point on $[a,c]$ at distance at most $M_1$ from $d$. Then by (H1) and definitions of $d,d'$: 
$$\grom{d}{c}{b} = \frac{1}{2} (\abs{b-d}+\abs{b-c}-\abs{c-d})=\frac{1}{2} (\abs{a-b}-\abs{a-d}+\abs{b-c}-\abs{c-d})= $$
$$= \frac{1}{2}(\abs{a-b}+\abs{b-c}-\abs{a-c})+ \frac{1}{2}(\abs{a-c} -\abs{a-d}-\abs{c-d}) \geq  $$
 $$\geq \grom{a}{c}{b} + \frac{1}{2}(\abs{a-c} -[\abs{a-d'}+\abs{d-d'}]-[\abs{d-d'}+\abs{c-d'}]) = \grom{a}{c}{b} -\abs{d-d'}.$$
 We obtain the first claim of part (ii) by using definition (H1) and the expression for $\grom{d}{c}{b}$ above 
 :
$$\abs{d-b}= \grom{d}{a}{b}\geq min \{\grom{d}{c}{b}, \grom {a}{c}{b}\}-\delta \geq   \grom{a}{c}{b}-\delta - \abs{d-d'},$$
and apply it to obtain the second claim:
$$\abs{a-c}=\abs{a-b}+\abs{b-c}-2 \grom {a}{c}{b} \geq \abs{a-b}+\abs{b-c}-2(\abs{d-b} +\abs{d-d'}+\delta).$$

(iii) Assume $d \notin B_{4\delta} ([a,c]),$ then $d \in B_{4\delta}([b,c])$ by (H2). We  apply   part (ii) to the points $b,a,c,d$ with $M_1=4\delta$ and obtain that $\abs{a-d} \geq \grom{b} {c} {a} -\delta-4\delta$. Contradiction. 

(iv) By definition (H3) we have that 
$$ \abs{d-c}+\abs{b-a} \leq max \{ \abs{a-c} +\abs{b-d},\abs{b-c} +\abs{a-d} \}+2\delta, \text{ hence }$$
$$\abs{d-c} \leq max \{ \abs{a-c} +(\abs{b-d}-\abs{b-a}),\abs{b-c} +(\abs{a-d}-\abs{b-a}) \}+2\delta \leq $$
$$ \leq max \{ \abs{a-c} ,\abs{b-c}  \}+2\delta +max\{\abs{b-d}-\abs{b-a}, \abs{a-d}-\abs{b-a} \}\leq$$ 
$$ \leq max \{ \abs{a-c} ,\abs{b-c}  \}+2\delta-D.\Box$$ 

\begin{lemm} \lab{lemma2*} Consider subgroups $\HH_1$ and $\HH_2$ in a finitely generated group $G$. If there exists $M \geq 0$ such that $\num{B_M(\HH_1) \cap B_M(\HH_2)}=\infty$ then  $\num{\HH_1\cap \HH_2}=\infty$.
\end{lemm}  
The lemma is equivalent to the statement: if $\num{\HH_1\cap \HH_2}<\infty$ then the set $B_M(\HH_1) \cap B_M(\HH_2)$ is finite for every non-negative $M$.

{\bf Proof} Assume that $\num{B_M(\HH_1) \cap B_M(\HH_2)}=\infty$ for some $M \geq 0$, then there exist infinite sequences of elements $\{h_{1i}\} \subset \HH_1$ and $\{h_{2i}\} \subset \HH_2$ such that $\abs{h_{1i}^{-1}h_{2i}}=\abs{h_{1i}-h_{2i} }\leq 2M$ for every $i \in \mathbb{N}$. We denote the element $h_{1i}^{-1}h_{2i}$ by $l_i$. 
Since $\abs{l_i} \leq 2M$ and the geometry of \cay \ is proper, there exists an element $l \in G$ and a subsequence $\{i_j \}$, $j \in \mathbb{N}$ such that $l_{i_{j_1}}=l_{i_{j_2}}=l$ in $G$ for any $j_1,j_2 \in \mathbb{N}$ and thus $h_{1s}^{-1} h_{2s}=h_{1k}^{-1}h_{2k}$ in $G$ for every $s,k \in \{i_j\}$. 
We obtained that $h_{1k}h_{1s}^{-1}=h_{2k}h_{2s}^{-1}$ belongs to  $ \HH_1 \cap \HH_2$ for every $ s,k \in \{i_j\}$. For every fixed $k$, $lim_{s\rightarrow \infty}\abs{h_{1k}h_{1s}^{-1}}\geq (lim_{s\rightarrow \infty} \abs{h_{1s}^{-1}})-\abs{h_{1k}}=\infty$ which imples that the intersection $\HH_1 \cap \HH_2$ does not belong to a ball $B_R$ for any $R \geq 0$ and thus is infinite. $\Box$ \medskip 

\begin{lemm} \lab{lemmaTFAE} Let $\HH$ be a $K$-quasiconvex subgroup in a $\delta$-hyperbolic group $G$ and let $E$ be an infinite elementary subgroup  in $G$. Then the following assertions are equivalent:

(i) for any number $M>0$ we have $E \not \subset B_M(\HH);$

(ii)$\num{E\cap \HH}<\infty $;

(iii) There exists $M>0$ (depending on $E$ and $\HH$ only) such that $\gro{x}{h}<M$ for any $x \in E$, $h \in \HH.$
\end{lemm}  
{\bf Proof}  We first show that (ii) implies (i). If the intersection $E \cap \HH$ is finite then by lemma \ref{lemma2*} we have  $\num{B_M(E) \cap B_M (\HH)}< \infty$  for any $M \geq 0$. In particular, $\num{E \cap B_M (\HH)}< \infty$ and hence $E \not \subset B_M(\HH)$ for any $M \geq 0$.

Now we show that (iii) implies (ii).  Let $x$ be an element of $E$ and $\abs{x}>2M$. We have 
$$\abs{x-h} =\abs{x} +\abs h -2\gro{x}{h} \geq \abs{x} -2M>0,$$
hence $x \notin \HH$ and so the intersection $E \cap \HH$ belongs to the ball $ B_{2M} (e)$ which is a finite set. 

It remains to show that (i) implies (iii). Since $E$ is infinite virtually cyclic we can choose an element $x$ in $E$ of infinite order and thus  $E$ is of finite index in $E(x)$. Hence 
\begin{equation} \lab{*M_0}
\text{there exists a constant } M_0 \text{ such that } E(x) \subset B_{M_0}(\langle x \rangle). 
\end{equation}
Let $K_x=K(\langle x \rangle)$ be a constant provided by lemma \ref{lemma 1.11}(iii). 

 Assume (iii) does not hold, i.e. for every  $M \geq 0$ there exist $y \in E$, $h \in \HH$, satisfying $\gro{y}{h}>M+M_0$. Then by (\ref{*M_0}), $y=x^{t}a$ for some $a \in E, \abs{a} \leq M_0$, $t \neq 0$ and 
  \begin{equation} \lab{TFAE1}
 \gro {x^{t}}{h} \geq  \gro {x^{t}a}{h}- M_0> M+M_0-M_0=M.
\end{equation}
Hence for every $M>0$ there exists an integer $t$ and an element $h \in \HH$ such that $\gro {x^{t}}{h} >M.$ 

Now we fix an arbitrary $t$ and choose $M$ so that  $\abs{x^{t}} <M-K_x-5\delta$. We may assume without loss of generality that $t \geq 0$. Then by (\ref{TFAE1}), there exist $t' \geq t$ and $h \in \HH$ such that $\gro {x^{t'}}{h}>M$. By lemma \ref{lemma 1.11}(ii) vertices $x^m$  are within $K_x$-neighborhood of $[e,x^{t'}]$ for any 
$0 \leq m \leq t'$. In particular, there exists a vertex $b \in [e,x^{t'}]$ such that $\abs{x^t-b} \leq K_x$ and thus $\abs{b} \leq M-5\delta < \gro {x^{t'}}{h}-5\delta$. By lemma \ref{lemma1*}(iii), we have that $b \in B_{4\delta}([e,h])$ and, because $\HH$ is $K$-quasiconvex, $b \in B_{4\delta+K}(\HH)$. Finally, we get that  $x^t $ belongs to $ B_{4\delta+K+K_x}(\HH)$ for every $t$ contrary to (i). $\Box$ \medskip    

In lemma \ref{lemma3*} and theorem \ref{quasiconvex and isomorphic} we follow in part the line of argument from \cite{Arzh} (in particular we apply lemma 13\cite{Arzh}). 
\begin{lemm} \lab{lemma3*} Let $x$ be an element of infinite order in 
$G$ and choose a constant $M_1 \geq 0$.  Then there exist a natural number $m$ and a number $M_2\geq 0$  such that for any element  $h$  in $G$ satisfying conditions $\abs{h} <2M_1$ and $h \notin E(x)$ and any $\abs t, \abs s \geq m$ the following inequality holds:
$$\abs {x^t h x^s} \geq \abs{x^t} +\abs{h} +\abs{x^s}-M_2.$$
\end{lemm}  

{\bf Proof}  For a pair of integers $s,t$ we consider a closed path $p_1q_1p_2 q_2$ in \cay, where the path $p_1$ starts  from $e$ and $lab(p_1) = x^{-t}$,  the path $q_1$ is geodesic  and  ends at vertex $h x^s$, the path $p_2$ satisfies $lab(p_2)= x^{-s}$ and $q_2$ is geodesic with $lab(q_2)= h^{-1}$. We define phase  vertices  $a_i$  on $p_1$ and phase vertices  $b_j$  on $p_2^{-1}$ relative to the natural factorizations $x^{-t}$ and $x^{s}$ respectively ($i=0, ...,-t$ and $j=0,...,s$). 
\medskip

Step 1. We take  constants $\lambda,c,K_x=K(\langle x \rangle)$ provided by lemma \ref{lemma 1.11} for the cyclic group $\langle x \rangle$ and define 
$$C=max \{2K_x+\frac{1}{2}\abs x , K_x+2M_1\}+8\delta.$$
 Let us denote by $y_i$ a phase path connecting vertex $a_i$ with some phase vertex $b_j$ of $p_2$. 
Assume that $\abs{y_{i}} \leq C$ for some $i$.   We define subpaths $p_1', p_2'$ of paths $p_1,p_2$, where the path  $p'_1$ connects $a_0$ to $a_{i}$ and $ p_2'$ connects $b_{j}$ to $b_0$. Considering the closed path $p_1'y_ip_2'q_2$, we have  
$$\abs j \abs x \geq \abs{h-h x^j} \geq \abs{x^i} -\abs{h} -\abs{y_i} \geq \lambda \abs x \abs i -c-2M_1-C,$$
 which implies that    $\abs{j} \geq \frac{\lambda \abs i \abs{x}-c-C-2M_1}{\abs x} \geq \lambda \abs i-c_1$, where $c_1= \frac{c+C+2M_1}{\abs x}$. 

 \begin{figure} \caption {}   \label{gamma}     
 \def\svgwidth{100mm}    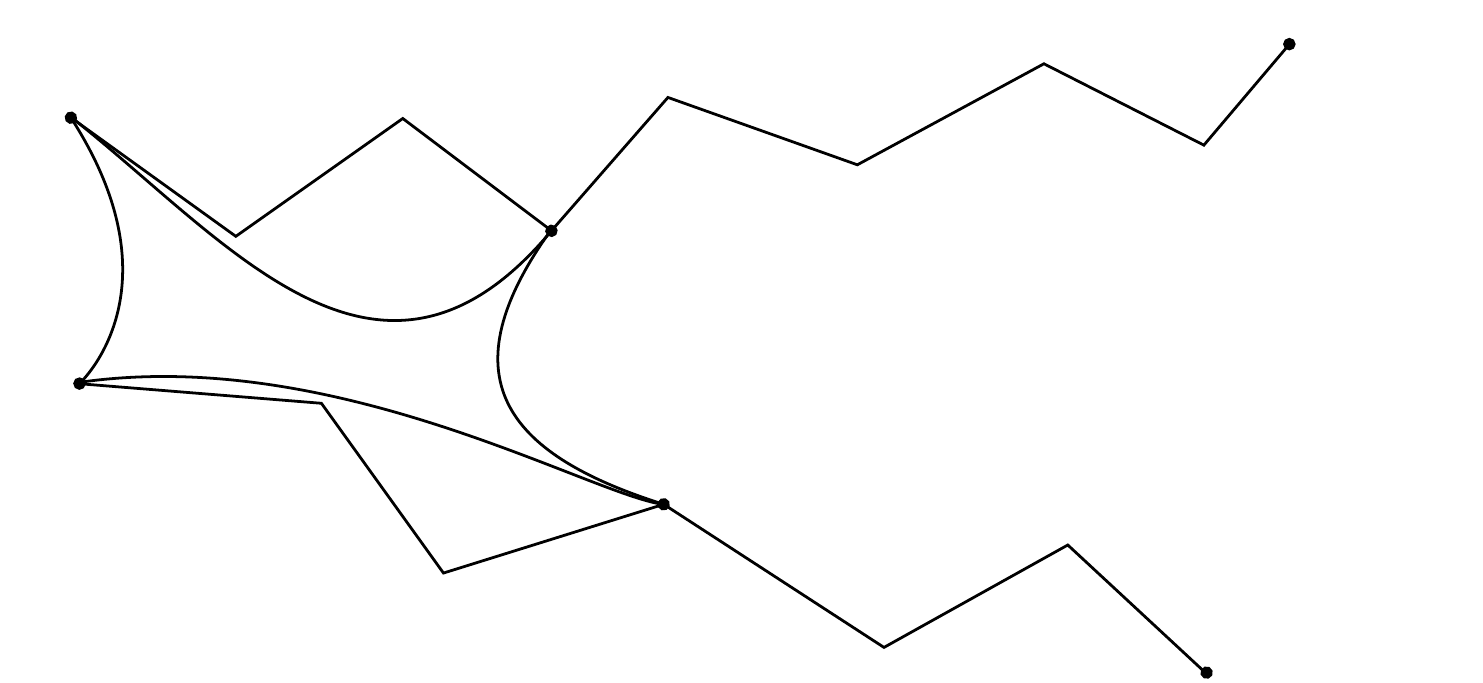 
 \end{figure} 
  
Since we have  fixed the constant $C$, we may apply lemma \ref{lemma 2.5} to the closed path $p_1'y_i p_2' q_2$ to obtain an integer $m_0$ such that if we choose a  number $i_0  $ satisfying $\lambda \abs{i_0}-c_1 \geq m_0$ and hence $\abs{i_0} \geq m_0$
 then there exists a phase path $y_{i'}$, $\abs {i'} \leq i_0$ such that $lab(y_{i'}) \in E(x)$. If the vertex $b_{j'}$ is the end vertex of $y_{i'}$, we get that $x^{-i'} lab(y_{i'})x^{-j'}h=e$ in $G$ and hence $h \in E(x)$, contradiction. We obtained that there exist $i_0$ depending on $x$ and $C$, such that 
\begin{equation} \label{cond i0}
\abs{y_{i_0}} > C.
\end{equation}
\medskip
Step 2. We show now that $a_{i_0} \in B_{8\delta+K_x}( q_1)$. By lemma \ref{lemma 1.11} $a_i \in B_{K_x}([e,x^{-t}])$ and using  twice the condition  (H2), we get that  $a_{i_0}$ belongs to $ B_{8\delta+K_x}(q_2 \cup[h, hx^s] \cup q_1)$. 

Clearly $a_{i_0} \notin B_{8\delta+K_x}(q_2)$: since $\abs{y_{i_0}}$ is minimal, i.e. $\abs{y_{i_0} } \leq \abs{a_{i_0}-b_{j'}}$ for every $j =0,...,s$, we get for $j'=0$ that    
$$\abs{y_{i_0} } \leq \abs{a_{i_0}-b_0} = \abs{a_{i_0}-h} \leq d(a_{i_0},q_2)+\abs{q_2} \leq 8\delta +K_x+\abs h <C$$
 contrary to (\ref{cond i0}).  

Similarly,  $a_{i_0} \notin B_{8\delta+K_x}([h, hx^s] )$. Otherwise we may consider a vertex $z$ on $[h,hx^s]$ at distance at most $8\delta+K_x$ from $a_{i_0}$ and choose a vertex $z'$ on $p_2$ at distance no more then $K_x$ from $z$. Finally, there exists a phase vertex $b_{j'}$ on $p_2$ such that $\abs{b_{j'}-z'} \leq \abs x /2$. Using the minimality of $\abs{y_{i_0}}$ we obtain the  estimate for the length of phase path: 
$$\abs {y_{i_0}} \leq \abs{a_{i_0}-z}+ \abs{z-z'} +\abs{z'-b_{j'}} \leq (K_x+8\delta)+K_x +\abs x /2 \leq C,$$
which again contradicts (\ref{cond i0}). The claim of Step 2 is proved. \medskip

Step 3. Let us choose some $\abs t , \abs s \geq \abs{i_0}$. We choose $z$ on $[e,x^{-t}]$ such that
\begin{equation} \lab{lemma3*1}
 \abs{a_{i_0}-z} \leq K_x.
 \end{equation}
 
   By Step 2 the vertex  $a_{i_0}$ is in the set $ B_{8\delta+K_x}( q_1)$ and hence $z \in B_{8\delta+2K_x}( q_1)$. Applying lemma \ref{lemma1*} (ii) to vertices $e,x^{-t}, hx^s, z$ we get 
   using (\ref{lemma1*2}), $\abs{h} <2M_1$: 
  $$\abs{q_1} \geq \abs{x^t}+\abs{hx^{s}}-2\abs{z}-2\delta-2(8\delta+2K_x)  \geq 
\abs{x^t}+ (\abs{h}+\abs{x^{s}}-4M_1)-2\abs{z}-2\delta-2(8\delta+2K_x).$$
Inequality  (\ref{lemma3*1}) implies that $\abs{x^{i_0} }+K_x \geq \abs z$ and we conclude that 
$$\abs{q_1} \geq \abs{x^t}+ \abs{h}+\abs{x^{s}}-2(\abs{x^{i_0}}+K_x)+2(9\delta+2K_x )=
\abs{x^t}+ \abs{h}+\abs{x^{s}}-(4M_1+2\abs{x^{i_0}}+6K_x+18\delta ).$$ 

It remains to define the constant $M_2$ (depending only on $\langle x \rangle$, $M_1$, $\HH$) to be $ 4M_1+2\abs{x^{i_0}}+4K_x+18\delta$ and define $m=\abs{i_0}$. $\Box$ \medskip
  
  \begin{lemm}(\cite{Arzh}, lemma 13) \lab{Arzh13} 
Let $n\geq 1$, $r\geq 48\delta$  and elements $h_i, g_i \in G$ ($1 \leq i \leq n$) satisfy :
\begin{equation} \lab{13*1}
\abs{g_i}>15r, \  (1 \leq i \leq n), \  \abs{h_1 g_1} \geq \abs{h_1} + \abs{g_1} - 2r,
  \end{equation}
\begin{equation} \lab{13*2}
\abs{g_{i-1}h_i g_i} \geq \abs{g_{i-1}}+\abs{h_i}+\abs{g_i}-2r (1<i \leq n ).
  \end{equation}

Then the following assertions are true:

(i) One has 
$$ \abs{h_1g_1 ... h_n g_n} \geq \abs{h_1g_1 ... h_{n-1}g_{n-1}} +\abs{h_n}+ \abs{g_n} -5r.$$
In particular one has (by induction) $h_1g_1 ... h_ng_n \neq e$ in $G$.

(ii) Let $p$ be a path in \cay \ labeled by  $h_1g_1 ... h_ng_n$. Then the path  $p$ and any geodesic $[p_-,p_+]$ are contained within $4r$-neighborhood of each other.
\end{lemm}

\begin{thm} \lab{quasiconvex and isomorphic} Let $G$ be a non-elementary hyperbolic group and $\HH$ be a $K$-quasiconvex subgroup of $G$. Consider an element $x$ in $G$ of infinite order such that $E(x)\cap \HH=\{e\}$. Then there exists  a number  $ r_0$ (depending on $\HH$ and $x$ only) such that 

(i) $\gro{x^s}{h}<\frac{r_0}{2}$ for any $h \in H$ and $ s \in \mathbb{Z}$ and 

(ii) for any $r \geq r_0$ there exists $t'>0$ such that for every $t \geq t'$ and   $g=x^t$   the subgroup $\HH_1=\langle g,\HH \rangle$ is $(4\delta+4r+max\{K,\abs g /2\})$-quasiconvex, of infinite index and canonically isomorphic to $\langle g \rangle*\HH$. Moreover, the inequalities of lemma \ref{Arzh13} hold for $r$; $g_i \in \ \langle g \rangle \backslash \{ e \}$ and $h_i \in \HH$.  
\end{thm}

{\bf Proof} 
By lemma \ref{lemmaTFAE}(iii) there exists $M>0$ such that $\gro{x^s}{h}<M$ for any $h \in H$ and $ s \in \mathbb{Z}$, hence
\begin{equation} \lab{eq2*}
\abs{hx^s}=\abs h +\abs {x^s}-2 \gro{x^s}{h} \geq \abs h +\abs {x^s}-2M.
\end{equation}

Now we consider an arbitrary element $ x^{m_1}hx^{m_2}$. If $\abs h > 2M+\delta$, then apply lemma \ref{lemma1*}(i) to vertices $e$, $x^{m_1}$,  $x^{m_1}h$, $x^{m_1}hx^{m_2}$ in \cay:  
\begin{equation} \lab{eq3*}
\abs{ x^{m_1}hx^{m_2} } \geq \abs{ x^{m_1}} + \abs h + \abs{ x^{m_2}}-4M-2\delta.
\end{equation}
 lemma \ref{lemma 1.11}(iv) provides a constant $M' \geq 0$ such that the following inequality  holds provided $m_1m_2 \geq 0$: 
\begin{equation} \lab{eq3*+}
\abs{ x^{m_1}x^{m_2} } \geq \abs{ x^{m_1}}  + \abs{ x^{m_2}}-2M'.
\end{equation}  
By lemma \ref{lemma3*}, there exist a natural number $m$ and a non-negative constant $M_2$ such that for any $\abs{m_1},\abs{m_2} \geq m$ the following inequality holds:
\begin{equation} \lab{eq4*}
\abs{ x^{m_1}hx^{m_2} } \geq \abs{ x^{m_1}} + \abs h + \abs{ x^{m_2}}-M_2 
\end{equation}
 for every  $h\neq 1, \abs{h} <2M+2\delta$. Now we choose
 \begin{equation} \lab{choose r}
  r_0=max\{2M+\delta, M_2/2, 48\delta, M' \}, 
\end{equation}  

and then for any $r \geq r_0$ we choose $t'$ satisfying $\abs {t'} \geq m$  so that the inequality
 \begin{equation} \lab{choose t}
    \abs{x^{t}} >15 r, \text { holds for every } t \abs t \geq \abs{t'}. 
\end{equation}

 For $g=x^t$   we consider an arbitrary element $h_1g^{s_1} ... h_n g^{s_n}$, where $n \geq 1$, $s_1 ... s_n \neq 0$ and every $h_i \in \HH \backslash \{e\}$ for $i=1, ..., n$.  We check the first condition of lemma \ref{Arzh13}, i.e. $\abs{g^s} >15 r $. For $s=1$ it is provided by (\ref{choose t}), if $s>1$ then:
$$\abs{g^s} = \abs{x^{st}} \geq \abs{x^{(s-1)t}}+\abs {x^t} -2M'> \abs{g^{s-1}}+13M'>15r.$$ 
The second  condition of lemma \ref{Arzh13} is  satisfied by (\ref{eq2*}) and the third because  (\ref{eq3*})--(\ref{eq4*}) hold.

We obtain that equality $h_1g^{s_1} ... h_n g^{s_n}=e$ in $ G $,  where $n \geq 1$, implies that either $s_1 ... s_n = 0$  or $h_i = e$  for some $ i=1, ..., n.$ Thus the group generated by $\HH$ and $g$ is isomorphic to $\HH*\langle g \rangle$. 

Consider an element $h=h_1g^{s_1} ... h_n g^{s_n}h_{n+1}$ in $G$ where $s_1 ... s_n \neq 0$ and $h_i\neq e$ for $i \leq n$ ($h_{n+1} $ can be the identity $e$). Define a  path $pp'$ in \cay \  with $p$ starting from  $e$ and label $h_1g^{s_1} ... h_ng^{s_n}$ in \cay \ and path $p'$ labeled by $h_{n+1}$. By (H2) and lemma \ref{Arzh13} (ii)  we have that $[e, p'_+] \subset B_{4\delta+4r}(p \cup p')$. In turn, every vertex $v$ of $p \cup p'$ is either within $K$--neighborhood of $\langle \HH,g \rangle$ (if $v$ is a vertex of a subpath labeled by $h_i$) or at most $\abs{g}/2$ away from $\langle\HH,g \rangle$ (if $v$ is a vertex of the subpath labeled by $g^{s_i}$). We conclude that every vertex $z \in [e, p'_+]$  is within $4\delta+4r + max \{ K,\abs g/2 \}$ from a vertex in $\HH_1=\langle\HH,g\rangle$. Hence $\HH_1$ is  ($4\delta+4r + max \{ K,\abs g/2 \}$)--quasiconvex. All conclusions of the theorem are checked for $\HH_1$ except the infiniteness of index. It only remains to observe that the subgroup $\HH_2=\langle\HH,g^2\rangle$ has infinite index in  $\langle\HH,g\rangle$ and hence in $G$. It satisfies the all of the conditions of the theorem and hence the conclusion holds for the same constant $r$ and  $g=x^{2t}$.$\Box$ \medskip
 
 Let us consider a path $p$ in \cay \ starting at vertex $a$ and ending with $b $ with label $h_1g^{s_1} ... h_n g^{s_{n}} h_{n}$, i.e. $ah_1g^{s_1} ... h_n g^{s_{n}} h_{n}=b$ in $G$ where $s_i \in \{\pm 1 \}$, $h_i \in \HH$ and if $h_i=e$ for $1 < i \leq n $ then $s_{i-1} s_i=1$. 
 We shall denote:
\begin{equation} \lab{a-i,b-i}
\begin{array} {cc} 
a_0=a, &  b_1=ah_1,    \\
   a_i=ah_1 ... h_ig^{s_i}  \text{ for } 1<i\leq n,  &  b_i=ah_1g^{s_1} ... h_i,\text{ for } 1<i\leq n.
 \end{array} 
 \end{equation}
     
 \begin{lemm} \label{aboutdistance} Assume that theorem \ref{quasiconvex and isomorphic} holds for $\HH$,$x$. Take some  constant $r$ and an element $g=x^t$ satisfying the same theorem. 
  In the notations (\ref{a-i,b-i})  we have that  $\grom{a} {b_{i+1}} {{a_i}} \leq r+\delta$              
  for any $i\geq 1$. 
  \end{lemm}
   {\bf Proof} The definition (H3) for $a,b_i,a_i,b_{i+1}$ reads: 
  \begin{equation} \lab{aboutdistance1}
  \abs{a_i-a}+\abs{g^{s_i}h_{i+1}} \leq max \{ \abs{b_i-a}+\abs{h_{i+1}}, \abs{b_{i+1}-a}+\abs{g^{s_i}}     \} +2\delta.
  \end{equation}
   
   The  theorem \ref{quasiconvex and isomorphic} permits us to apply the inequalities of lemma \ref{Arzh13}(i)  and the second condition in (\ref{13*1}) to the left side of (\ref{aboutdistance1}): 
 $$\abs{a_i-a}+\abs{g^{s_i}h_{i+1}} \geq (\abs{a_{i-1}-a}+\abs{h_i}+\abs{g^{s_i}}-5r)+ (\abs{g^{s_i}} +\abs{h_{i+1}}-2r),$$ 
 applying  the first inequality in (\ref{13*1}) we obtain
 $$\abs{a_i-a}+\abs{g^{s_i}h_{i+1}} \geq (\abs{a_{i-1}-a}+\abs{h_i}) +\abs{h_{i+1}}+(2\abs{g^{s_i}}-7r) > \abs{b_i-a}+\abs{h_{i+1}} +23r.$$ 
  By the conditions on $r$ in theorem \ref{quasiconvex and isomorphic}:
  \begin{equation} \lab{aboutdistance2}
 \abs{a_i-a}+\abs{g^{s_i}h_{i+1}} > \abs{b_i-a}+\abs{h_{i+1}} +23r>\abs{b_i-a}+\abs{h_{i+1}} +2\delta. 
  \end{equation}
   
  Hence we may rewrite (\ref{aboutdistance2}) as    $\abs{a_i-a}+\abs{g^{s_i}h_{i+1}} \leq \abs{b_{i+1}-a}+\abs{g^{s_i}} +2\delta$ and thus (using the second inequality of (\ref{13*1}) again):
 $$\abs{b_{i+1}-a}\geq \abs{a_i-a}+\abs{g^{s_i}h_{i+1}} -\abs{g^{s_i}} -2\delta \geq \abs{a_i-a}+\abs{g^{s_i}} +\abs{h_{i+1}}-2r -\abs{g^{s_i}} -2\delta = $$ 
 $ =\abs{a_i-a} +\abs{h_{i+1}}  -2\delta -2r$, which by (H1) implies that $\grom{a} {b_{i+1}} {{a_i}} \leq r+\delta. \Box$ \medskip
 
  \begin{lemm} \lab{path and ball}  Assume that theorem \ref{quasiconvex and isomorphic} holds for $\HH$, $x$. Take some  constant $r$ and an element $g=x^t$ satisfying this theorem.
  In the conventions above (see (\ref{a-i,b-i})), assume that 
 
 (i) $a_i \in B_{R}$ for some $i>1$. Then $a_1$ belongs to $B_{R+2\delta-2r}$.  

 (ii) $b_{i+1} \in B_{R}$ for some $i >1$.  Then $a_1$ belongs to $B_{R}$.
 \end{lemm}
 
 {\bf Proof} (i) By lemma \ref{Arzh13}(ii) there exists $b' \in [a,a_i]$ such that 
 \begin{equation} \lab{path and ball01}
 \abs{b'-a_1} \leq 4r.
 \end{equation}
  Using the inequality (\ref{13*1}) of the same lemma, we have 
 
\begin{equation} \lab{path and ball02}
\abs{b'-a} \geq \abs{a-a_1} - \abs{b'-a_1} \geq  \abs{g} +\abs {h_1}-2r -4r \geq 9r.
\end{equation}

  Similarly, we may inductively apply lemma \ref{Arzh13}(i) to the subpath of $p$ connecting $a_j, a_i$ for $j<i$:
    \begin{equation} \lab{path and ball1}
  \abs{a_j-a_{i}} \geq \abs{a_j-a_{i-1}} +  \abs{h_i} + \abs{g^{s_i}}-5r  > \abs{a_j-a_{i-1}} + 15r-5r \geq 10r(i-j),
   \end{equation}
   and apply it  in order to estimate (for $i>1$): 
  $$\abs{b'-a_i} \geq \abs{a_1-a_{i}}-\abs{a_1-b'} \geq 10r(i-1)-4r \geq 6r.$$
  
   The inequalities (\ref{path and ball01}) and (\ref{path and ball02}) allow  to apply lemma \ref{lemma1*}(iv) to $a,$ $a_i$, $e$, $b'$ (with $D=6r$) and get that 
    $\abs{b'-e} \leq max\{ \abs{a-e},\abs{a_i-e}\}+2\delta-6r .$ Since $a,a_i \in B_R$ we have that $\abs{b'-e} \leq R+2\delta-6r$. We use the previous inequality together with (\ref{path and ball01}) to conclude that:
   $\abs{a_1-e} \leq \abs{a_1- b'} + \abs{b'-e}  \leq 4r +R-6r+2\delta=R+2\delta-2r.$ \medskip
     
 (ii) The inequality (\ref{path and ball1}) we have that $\abs{a_i-a} \geq 10r >r+6\delta+1$; on the other hand lemma \ref{aboutdistance} implies that 
 $r+6\delta+1 \geq \grom{a} {b_{i+1}} {{a_i}}+5\delta+1$.  Thus we can choose a vertex $d$ on a geodesic $[a,a_i]$ satisfying inequalities:
  $$\grom{a}{b_{i+1}}{{a_i}}+5\delta+1 \geq \abs{d-a_i} \geq  \grom{a} {b_{i+1}} {{a_i}}+5\delta.$$
 
   Then, by lemma \ref{lemma1*}(iii), $d$ belongs to $B_{4\delta} ([a,b_{i+1}])$ and using  lemma \ref{aboutdistance} we get 
  $$ d(a_i, [a,b_{i+1}]) \leq \abs{d-b}+4\delta \leq \grom{a}{b_{i+1}}{{a_i}}+5\delta+1+4\delta \leq r+10\delta+1.$$
 By (H2), segment $[a,b_{i+1}]$ belongs to the $4\delta$--neighborhood of  union $[e,a] \cup [e,b_{i+1}]$  which is a subset of $B_R$ because  $a,b_{i+1} \in B_{R}$. Hence 
 $$\abs{a_i-e} \leq d(a_i, [a,b_{i+1}])+4\delta+R \leq R+r+14\delta+1$$
    and by part (i) of this lemma we conclude that $a_1$ belongs to $B_{R- r+16\delta+1} \subset B_{R}. \Box$ \medskip 
      
\begin{lemm} \lab{ballprop} Let $\HH$ be a $K$-quasiconvex subgroup in a hyperbolic group $G$. Assume that $a \in B_R$ and $ah \notin B_R$ for some $h \in \HH$. Then either $\gro{a^{-1}}{h} \leq 13\delta+K$ or there exists $b_1 \in aH \cap B_R$ such that $b_1h_1=ah$ and $\abs{h_1} < \abs{h}$ for some $h_1 \in \HH$.
\end{lemm} 
{\bf Proof} Assume that    $\gro{a^{-1}}{h} > 13\delta+K$. We choose a vertex $d$ on the segment $[a,ah]$ such that $\abs{d-a}=K+8\delta$. By lemma \ref{lemma1*}(iii), $d \in B_{4\delta}([e,a])$ and we can choose $d' \in [e,a]$ to satisfy the inequality $\abs{d-d'} \leq 4\delta$. Then we have  $$\abs{d-e} \leq \abs{e-d'}+\abs{d-d'} \leq (\abs{e-a}-\abs{a-d'})+4\delta \leq R-\abs{a-d'} +4\delta \leq $$
$$ \leq  R-\abs{a-d} +4\delta  +4\delta \leq R-K.$$

By quasiconvexity of $\HH$, there exists $b_1 \in a\HH$, $\abs{b_1-d}\leq K$ and hence $b_1 \in B_R$. By the choice of $b_1$ we have that  $b_1 ^{-1} ah=h_1 \in \HH$ and 
$$\abs{b_1-ah} \leq \abs{b_1-d}+\abs{d-ah} \leq \abs{b_1-d}+(\abs{a-ah}-\abs{d-a}) \leq K+ (\abs h-K-8\delta)<\abs h. \Box$$   

  \begin{lemm} \lab{ballprop2}  Assume that theorem \ref{quasiconvex and isomorphic} holds for $\HH$, $x$. Take some  constant $r$ and an element $g=x^t$ satisfying this theorem. We adopt notations (\ref{a-i,b-i}) and let $a, b$ be vertices in $ B_R$ and $ah_1g^sh_2=b $ in $G$ for some $h_1,h_2 \in \HH$, $s \in \{\pm 1 \}$. Assume furthermore that  $\gro{a^{-1}}{h_1} \leq 13\delta+K$ and that $b_1 \notin B_R$. Then 
$$\abs{h_1} \leq K+\frac{r_0}{2}+15\delta.$$
 \end{lemm}     
{\bf Proof} Definition (H1) and theorem \ref{quasiconvex and isomorphic} yield:
\begin{equation} \lab{ballprop2*}
\frac{r_0}{2} \geq \grom{a_1}{a}{{b_1}} \geq min  \{ \grom{a}{b}{{b_1}}, \grom{a_1}{b}{{b_1}} \}-\delta \geq  min\{ \grom{e}{a}{{b_1}}, \grom{e}{b}{{b_1}}, \grom{a_1}{b}{{b_1}} \} -2\delta. 
\end{equation}     
     
Consider the last two Gromov products on the right-hand side of (\ref{ballprop2*}). We have: 
$$\grom{e}{b}{{b_1}}=\frac{1}{2}(\abs{{b_1}}+\abs{b-b_1} - \abs b)= \frac{1}{2}(\abs{{b_1}} - \abs b) +\frac{1}{2}\abs{b-b_1}, $$
by the conditions of this lemma $\abs{b_1} > R \geq \abs b$ and using theorem \ref{quasiconvex and isomorphic} we conclude 
$$\grom{e}{b}{{b_1}} \geq 0+\frac{1}{2}(\abs g+ \abs {h_2}-r_0) \geq \frac{1}{2} \abs g -\frac{r_0}{2} >7r \geq 7r_0.$$  
Similarly,
$$\grom{a_1}{b}{{b_1}} = \frac{1}{2}(\abs{g}+\abs{b-b_1} - \abs {h_2}) \geq  \frac{1}{2}(\abs{g} +( \abs g+\abs{h_2}-r_0) -\abs{h_2}) \geq \abs g-\frac{r_0}{2} \geq 14\frac{1}{2} r_0 .$$ 

Now we may rewrite (\ref{ballprop2*}) as $\frac{r_0}{2} \geq  \grom{e}{a}{{b_1}} -2\delta.$ Note that $\grom{e}{b_1}{a}=\gro{a^{-1}}{h_1} \leq 13\delta+K$, thus by definition of the Gromov product (\ref{Gromprod}): 
$$\abs{h_1} =   \grom{e}{a}{{b_1}}+ \grom{{b_1}}{e}{a} \leq (\frac{r_0}{2}+2\delta)+(K+13\delta).\Box$$   

In order to estimate the number of $\HH$-cosets in $B_R$ from below, we define $M_R=\{a\HH \vert \  a\HH \cap B_R \neq \emptyset  \}$ and $Q_R=\{a\HH \in M_{R} \vert \exists b \in B_r, \ b\HH \neq a\HH \ \& \ b \langle \HH,g \rangle = a \langle \HH,g \rangle  \}$.

\begin{lemm} \lab{estimateBR} Let $\HH$ be a free $K$--quasiconvex subgroup in $G$. Then for any $k \in \mathbb{N}$ and any $x \in G$ of infinite order either:

(i) there exists $t \neq 0$ such that $x^t \in \HH$, or

(ii)there exists $t$ such that for $g=x^t$ the group $\langle \HH,g \rangle$ is  quasiconvex and canonically isomorphic to the free product 
$\HH* \langle g \rangle$. Moreover $\frac{\num{Q_R}}{\num{B_R}} \leq \frac{1}{2^k}$ for any $R >0$. 
\end{lemm}
{\bf Proof} Assume that (i) does not hold, so $x^t \in \HH$ implies $t =0$. By lemma \ref{lemma2*} we have that for every $M \geq 0$ the number of vertices in $B_M(\langle x \rangle) \cap B_M(\HH)$ is finite.  There exists $M_0 \geq 0$ such that  $E(x)$ is in $M_0$--neighborhood of $\langle x \rangle$, hence  $B_M(E(x)) \cap B_M(\HH) \subset B_{M+M_0}(\langle x \rangle) \cap B_{M+M_0} (\HH)$ and hence $\num{B_M(E(x)) \cap B_M(\HH)}$ is finite thus  $\num{E(x) \cap \HH}< \infty$. Since $\HH$ is free, the last inequality means that $E(x) \cap \HH=\{e \}.$

Using corollary \ref{amenable}, we choose $c$ so that that $2^{k+1} (\num{B_{K+r_0+15\delta}}^2 \num{B_{R-c}}) \leq \num{B_R}$ and $r$ according to the theorem \ref{quasiconvex and isomorphic}  and satisfying:
\begin{equation} \lab{estimateBR for r}
r \geq max\{c/7,2(K+\frac{r_0}{2}+17\delta)\}.
\end{equation}

We choose $t$ according to theorem \ref{quasiconvex and isomorphic}. 

  Let $a\HH$ belong to $Q_R$, then there exist $b \in B_R$, $b \notin a\HH$, and  elements $h_i \in \HH$ ($i=1,...,k$) such that $ah_1g^{s_1} ... g^{s_k} h_{k+1}=b$ in $G$. By lemma \ref{path and ball} we have that either $b_1=a h_1$ or $a_1=a h_1g^{s_1}$ or $ah_1g^{s_1}h_2$ belongs to $B_R$. Hence we can assume that $b=a h_1g^{s_1}h_2$, where $h_1,h_2 \in \HH, \ s \in \{\pm 1 \}$. Clearly $a \HH\neq b\HH$, otherwise $ah_1g^{s_1}h_2=ah$ and $g^{s_1}=x^{ts_1} \in \HH$, contradiction.

We may also assume that $a, h_1$ are chosen so that $\abs{h_1}$ is minimal with respect to all factorizations $a'h'=ah_1$ in $G$ where $a' \in a\HH \cap B_R$.  Similarly, we may assume that $b,h_2$ are chosen so that for any $b' \in b\HH \cap B_R$ and $h' \in \HH$  the equality 
$b' {h'}^{-1}=b h_2^{-1}$ implies that $\abs{h'} \geq \abs{h}$. According to the choice we made, if $h_1\neq e$ ($h_2\neq e$) then $b_1=ah_1 \notin B_R$ (respectively $a_1=ah_1g^s \notin B_R$).  Now we are in position to apply lemma \ref{ballprop} to the pairs $a,ah_1$ and 
$b, bh_2^{-1}$, which provides that 
$\gro{a^{-1}} {h_1} \leq K+13\delta,\  \gro{b^{-1}} {h_2^{-1}} \leq K+13\delta,$
and then, by lemma \ref{ballprop2}, we conclude that 
\begin{equation} \lab{estimateBR2}
\abs{h_i} \leq K+\frac{r_0}{2}+15\delta,  \text{ for } i=1,2. 
\end{equation}
If $h_i=e$ for $i=1$ or $2$ then the corresponding inequality in (\ref{estimateBR2}) holds trivially.

 We have that $b_1, a_1 \in B_{R+K+\frac{r_0}{2}+15\delta}$. Since $\abs{g}>15r$, we can fix a factorization $g=g_1g_2$ in $G$ such that $\abs{g_1}+\abs{g_2}=\abs{g}$ and $\abs{g_1},\abs{g_2} \geq \frac{15r}{2}.$ Let  $b'=ah_1g_1$ if $s=1$ and $b'=ah_1g_2^{-1}$ if $s=-1$, we will call $b'$ a middle point of the path $p$ starting at $a$ with label  $h_1g^sh_2$. 

Applying lemma \ref{lemma1*}(iv) to vertices $b_1,a_1,b'$,  we obtain that $\abs{b'-e} \leq (R+K+\frac{r_0}{2}+15\delta) +2\delta -\frac{15r}{2}$. As we choose $r$ according to (\ref{estimateBR for r}) we obtain:
\begin{equation} \lab{estimateBR3}
b' \in B_{R-7r}. 
\end{equation}

We have obtained that if  a coset $a'\HH$ belongs to $Q_R$ then there exist $a \in a'\HH \cap B_R$, $b \in B_R$, $h_1, h_2 \in \HH$ and $s \in \{ \pm1\}$ such that the equation $ah_1g^sh_2=b$ holds in $G$ together with conditions (\ref{estimateBR2}) and (\ref{estimateBR3}). Hence the number of elements in $Q_R$ is not greater then the number of paths with label $h_1g^sh_2$ in \cay \  such that the middle point $b'$ of each path satisfies (\ref{estimateBR3}):
\begin{equation} \label{estimateBR4}
\num{Q_R} \leq  \# \{ \text{of  } h_1g^sh_2 \text { satisfying (\ref{estimateBR2})}\} \times \num{B_{R-7r}} \leq 2\num{B_{K+\frac{r_0}{2}+15\delta}}^2 \times \num{B_{R-7r}}. 
\end{equation}

Due to our choice of $r$ in (\ref{estimateBR for r}) we finally get 
$$\num{Q_R} \leq  2\num{B_{K+\frac{r_0}{2}+15\delta}}^2 \num{B_{R-7r}} \leq \frac{1}{2^k}\num{B_R}. \Box$$  

 \begin{rem} \lab{existance} Let 
 $\HH$ be an infinite quasiconvex subgroup of $G$ of infinite index. Then there exists an element $x \in G$  of infinite order  such that it is non-commensurable with any element of $\HH$.
 \end{rem}  
In particular, this remark implies that no infinite index subgroup satisfying the Burnside condition in a non-elementary hyperbolic group is quasiconvex. 

{\bf Proof} 
By proposition \ref{doublecoset}, for every $N>0$ there exists a double coset $\HH g \HH$ which has no representative of length shorter then  $N$. Choose an element $g$ for some $N>2K+2\delta$. Then by lemma \ref{width} the intersection $\HH \cap g^{-1} \HH g$ is finite. The subgroup $\HH$ contains an element $h$ of infinite order because (by Lemma 18, \cite{IvOl}) every torsion subgroup in a hyperbolic group is finite. The intersection  $\langle g^{-1}hg \rangle \cap \HH \subset  g^{-1} \HH g\cap \HH$  is finite and hence  $x=g^{-1}hg$ is non-commensurable with any element of $\HH$.$\Box$ \medskip 

\begin{thm} \lab{Burnsidecoset} For every non-elementary \hyp \ group $G$ and any $0<q<1$ there exists a free  subgroup $H$ satisfying the Burnside condition and  such that 
$\frac{ \# \{  aH \vert \  aH  \cap B_R \neq \emptyset \} }{\num{B_R}} \geq q$.  
\end{thm}
{\bf Proof} We choose a sequence $\{k_i \}_{i \in \mathbb{N} }$ such that 
\begin{equation} \label{series}
\Sigma_{i=1}^\infty \frac{1}{2^{k_i}} <1-q.
\end{equation}

Let $\{x_j\}$, $j \in \mathbb{N}$ be a list of all elements of infinite order in $G$. We fix notations $\HH_i=\langle x_1^{t_1}, ... , x_i^{t_i} \rangle$ for some positive numbers  $t_i \in \mathbb {N}$ which we will determine later. We define $H=\cup_{i=1}^\infty \HH_i $, it clearly satisfies the Burnside condition. Then we denote 
$M_R^i=\{a \HH_i \vert \ a \HH_i \cap B_R  \neq \emptyset \}$ and $Q_R^i= \{a \HH_i \vert \  \exists b \in B_R \text{ such that } a\HH_i \neq b\HH_i \ \& \ a\HH_{i+1} = b\HH_{i+1} \}$. 

We set $\HH_0=\{ e \}$ and thus $M_R^0=B_R$. lemma \ref{estimateBR} (applied to $\HH_0$, $x_1$ and $k_1$) provides that  there exists $t_1>0$ such that $\HH_1=\langle g_1 \rangle$ (where $g_1=x_1^{t_1}$) is cyclic, quasiconvex and $ \frac {\num{Q_R^0}} {B_R} \leq \frac{1}{2^{k_1}}$ for any $R > 0$. It provides the  following estimate for $M_R^1$:
$$\num{M_R^1}\geq \num{M_R^0}-\num{Q_R^0} \geq (1-\frac{1}{2^{k_1}}) \num{B_R}.$$

Now we assume by induction that a free quasiconvex subgroup $\HH_i=\langle x_1^{t_1}, ... , x_i^{t_i}  \rangle$ has been constructed by repeated application of lemma \ref{estimateBR} and $M_R^{i}$ satisfies inequality 
\begin{equation} \lab{induction}
\num{M_R^{i}}\geq  (1-\frac{1}{2^{k_1}}-\frac{1}{2^{k_2}}- ... -\frac{1}{2^{k_i}}) \num{B_R}.
\end{equation}
 If $\langle x_{i+1}\rangle \subset B_M(\HH_i)$ for some $M \geq 0$ then by lemma \ref{lemmaTFAE} there exists $t_{i+1}>0$ such that $x_{i+1}^{t_{i+1}} \in \HH_i$ and we can set $\HH_{i+1}=\HH_i$, finishing the induction step ($M_R^{i+1}=M_R^{i}$).

Assume now that $ \langle x_{i+1} \rangle \not \subset B_M(\HH_i)$ for any non-negative $M$,   then by lemma \ref{lemmaTFAE} we have  $\num{E(x) \cap \HH_i}< \infty$ and hence (because $\HH_i$ is free) $E(x) \cap \HH_i=\{e \}$. We choose $t_{i+1}$ appplying  lemma \ref{estimateBR} to $\HH_i,x_{i+1},k_{i+1}$ and using the induction assumption (\ref{induction}):
 $$ \num{M_R^{i+1}}\geq \num{M_R^i}-\num{Q_R^i} \geq (1-\frac{1}{2^{k_1}}-\frac{1}{2^{k_2}}- ... -\frac{1}{2^{k_i}}) \num{B_R}- \num{Q_R^i} \geq$$
 $$\geq (1-\frac{1}{2^{k_1}}-\frac{1}{2^{k_2}}- ... -\frac{1}{2^{k_{i+1}}}) \num{B_R},$$

and by (\ref{series}):
  $$\num{M_R} \geq (1-\Sigma_{i=1}^\infty \frac{1}{2^{k_i}})\num{B_R} >q\num{B_R}. \Box$$
  
{\bf Proof of theorem \ref{maxgrowthburnside}} By the remark \ref{symmetry}, the number of left cosets intersecting $B_n$ is equal to the number $f_{H\backslash G}(n)$ of right ones. We can now fix some $0<q<1$ and using theorem \ref{Burnsidecoset} find a group $H$ such that $f_{H\backslash G}(n) \geq q f(n).$ Thus (by remark \ref{morethenexponent}) the growth of action of $G$ on $H\backslash G$ is maximal. $\Box$

\section{Proof of theorem \ref{generalizedA3} and corollary \ref{generalizedA2}}
\begin{lemm} \lab{generalizedA1} Let $G$ be a non-elementary hyperbolic group and $\HH$ be a quasiconvex subgroup in $G$ such that $E(G) \cap \HH=\{e \}$, then there exists $g \in G$ of infinite order such that $E(g)=\langle g \rangle \times E(G)$ and $\langle \HH ,g \rangle $ is quasiconvex of infinite index and is canonically isomorphic to $ \HH *   \langle g \rangle.$
 \end{lemm}  

{\bf Proof} By lemma \ref{3.4,3.8} there exists $y \in G$ of infinite order such that $E(y)= \langle y \rangle\times E(G).$ We have either 

(I) $\num {E( y ) \cap \HH } <\infty,$ 

or (II) $\num{E( y ) \cap \HH } $ is infinite. 

Take an  element $y^k a \in E(g)$ and assume $y^k a \in \HH$ for some $k \in \mathbb Z$ and $a \in E(G)$, then  we have $(y^k a)^n=y^{kn} a^n$ for every $n$ (because $a$ commutes with $y$). Note that for a non-zero $k$ the equality  $y^{kn_1} a^{n_1}=y^{kn_2} a^{n_2}$ holds in $G$ if and only if $n_1=n_2.$

Hence in   case (I) we have that $k=0$ and thus $E( y ) \cap \HH \subset E(G) \cap \HH=\{e \}$ and we may apply theorem  \ref{quasiconvex and isomorphic} to find $t>0$ and obtain the canonical isomorphism  $ \langle \HH, y^t \rangle \cong \HH* \langle y^t \rangle.$

Thus we only need to consider case (II) when $y^k a \in \HH$ for some non-zero $k$. Replacing $y$ with $y^ka $ we may assume that $y$ is in $\HH$. By remark \ref{existance} there exists an element $x$ of infinite order such that $x$ is non-commensurable with any element of $\HH.$ Replacing $x$ with its non-zero power if necessary, we may assume that $x$ commutes with $E(G)$.
 We define a  subgroup $\HH_1=\langle y \rangle$, which is quasiconvex by lemma \ref{lemma 1.11}(ii). Since $\HH_1 \cap E(x)= \{ e \}$, there exists a constant $r_0 \geq 0$ such that $\gro {x^t} {y^s} < \frac{r_0}{2}$ for all $t,s \in \mathbb Z$ by theorem \ref{quasiconvex and isomorphic}.
By part (ii) of the same theorem and  $r=r_0$, there exists $t'>0$ such that $\langle x^{t'},  y \rangle \cong \langle x^{t'} \rangle *\langle y \rangle. $ We denote $x_1=x^{t'}$ so  the subgroup $\langle x_1,y \rangle$ is free quasiconvex and inequalities of lemma \ref{Arzh13} hold for $\HH_1,x_1$ and $r=r_0$. In particular, for every reduced word $w$ in $\langle x_1,y \rangle$, the corresponding path in \cay \ with label $w$ is within $4r_0$-neighborhood of a geodesic connecting its ends. 

By lemma \ref{lemmaTFAE} there exists $M \geq 0$ such that $  \gro {x^s} {h} <M$ for every $s \in \mathbb Z, \ h \in \HH$. 
Choose $t>0$ such that $\abs{x_1^s} > 4r_0+K+2M$ for every $\abs s \geq t$ and denote $x_2=x_1^t$.  We have that for any non-zero $m$:
\begin{equation} \lab{generalizedA1eq1} 
d (x_2^m, h) \geq \abs{x_2^{m}}+\abs h -2M \geq \abs{x_2^{m}}-2M > 4r_0+K.
\end{equation}
Let element $w=y^{s_0} x_2 ^{t_1} y ^{s_1} ... x_2^{t_n} y^{s_n}$   satisfy $s_i,t_i, t_n \neq 0$ for $i=1,...,n-1$ and assume that $w\in \HH$. Then every phase vertex of $w$ is in $(4r_0+K)$-neighborhood of $\HH$, which contradicts  inequality (\ref{generalizedA1eq1}) because 
 $ d(y^{s_0} x_2^{t_1}, y^{s_0})=d(x_2^{t_1}, e)>4r_0+K $. We conclude that an element $w$ of the free group $\HH_2=\langle x_2, y \rangle$ is commensurable with an element of $\HH$ if and only if $w=y^s$ for some integer $s$. 

Now we consider an element $y^kx_2^k$   for sufficiently large $k$ and will show that the group $E( y^kx_2^k)$ is equal to $\langle y^kx_2^k\rangle \times E(G)$. 
Let $z$ be an element of $E( y^kx_2^k)$, i.e. the equality   $z(y^kx_2^k)^m z^{-1}=(y^kx_2^k)^ { m'}$ holds in $G$ for some $m=\pm m'\neq 0$. 
We choose a constant
$$M_0 >2\abs z +8r_0+26\delta +(3k+1)(max\{\abs y, \abs{x_2} \})$$
 and a natural number $s$ divisible by $m$ such that $\abs{(y^kx_2^k)^s} \geq M_0$. We consider a closed path $p_1q_1p_2q_2$ in \cay \ such that $lab( p_1)=lab (p_2)=z$, $lab(q_1)= (y^k x_2^k)^s,$ and $lab(q_2^{-1})= (y^k x_2^k)^{s'},$  where $s'= \pm s$.  For convenience we denote the initial vertices of $p_1,q_1,p_2,q_2$ by $a,b,c,d$ respectively.
\begin{figure}   \caption {}     \label{5}   \def\svgwidth{120mm}  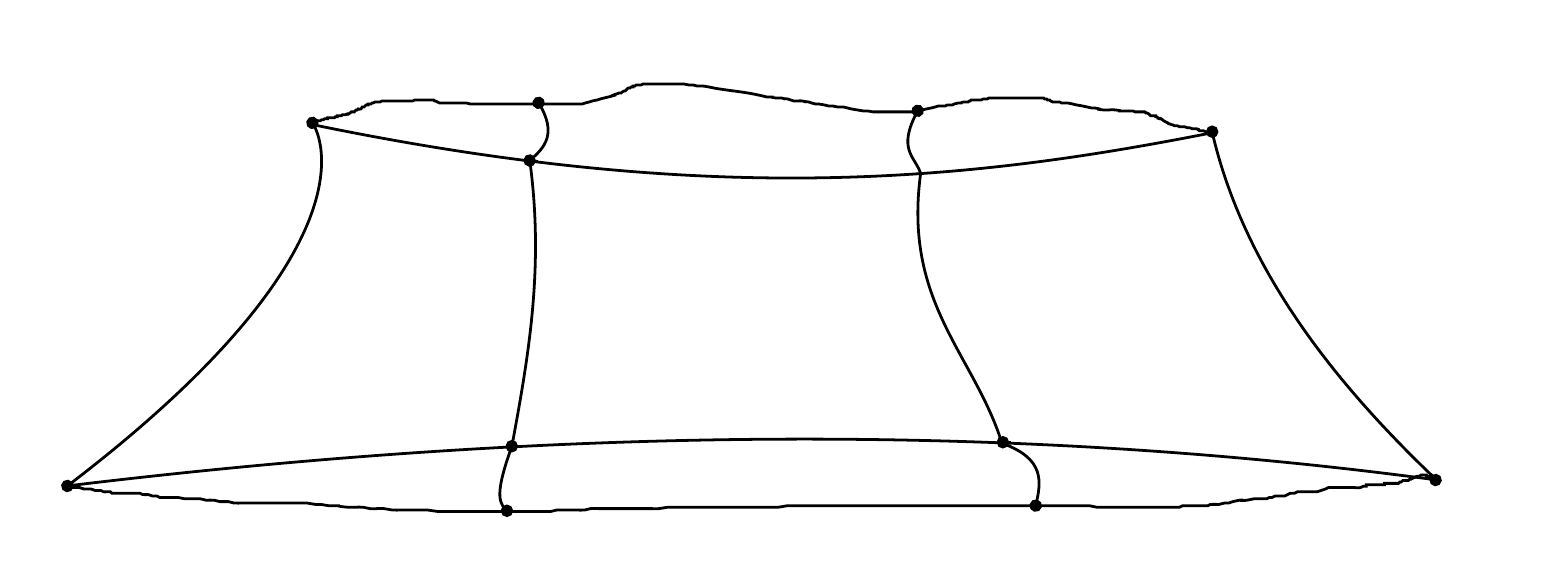  \end{figure} 
We choose a vertex $ \ov u$ on $[b,c]$  at distance $\abs z +5\delta$ in \cay \ from the vertex $b$. Then, using (H2), $\ov u$ is in $4\delta$-neighborhood of some $u_1$ on  $[a,c] \cup [a,b]$ and, by the choice of $\ov u$, is actually on $[a,c]$. Using (H2) again, and taking into account that 
$$\abs {\ov {u_1} -c} \geq \abs {b-c} -\abs z -5\delta-4\delta \geq M_0 -\abs z -9\delta>\abs{z}+ 5\delta$$
we obtain that there exists $\ov u '$ on $[a,d]$ satisfying $\abs{\ov u'-  u_1 } \leq  4\delta$. Hence 
\begin{equation} \lab{generalizedA1eq2} 
\abs{\ov u'- \ov u } \leq 8 \delta.
\end{equation}
 Similarly we can choose $ \ov v$ on $[b,c]$  at distance $\abs z +5\delta$ from the vertex $c$ and a vertex $\ov v'$ on  $[a,d]$ such that 
\begin{equation} \lab{generalizedA1eq3} 
\abs{\ov v'- \ov v } \leq 8\delta.
\end{equation} 
Since $q_1$ and $[b,c]$ are within $4r_0$-neighborhood of each other, we find phase vertices $u,v$ on $q_1$ relative to the factorization $y^k x_2^k...y^k x_2^k$ of $lab(q_1)$ such that 
\begin{equation} \lab{generalizedA1eq4} 
\abs{u- \ov u},\abs{v- \ov v}  \leq 4r_0+\frac{1}{2} max\{\abs y , \abs{x_2} \}.
\end{equation} 
Similarly we find phase vertices $u',v'$ on $q_2$ such that 
\begin{equation} \lab{generalizedA1eq5}
\abs{u'- \ov u'},\abs{v'- \ov v'}  \leq 4r_0+\frac{1}{2} max\{\abs y , \abs{x_2} \}.
\end{equation} 

Now we consider a closed path $p_1'q_1'p_2'q_2'$, where $q_1' ,q_2'$ are subpaths of $q_1$ and $q_2$ respectively and $p_1'=[u',u],$ $p_2'=[v',v]$. According to inequalities (\ref{generalizedA1eq2})--(\ref{generalizedA1eq5}) above:
\begin{equation}\lab{generalizedA1eq6}
\abs{p_i} \leq 8r_0+ max\{\abs y , \abs{x_2} \}+ 8\delta, \text{ where }i=1,2.
\end{equation}
Note that 
$$\abs{q_1'} = \abs{u-v} \geq \abs{\ov u- \ov v}- \abs{u- \ov u}-\abs{v- \ov v} \geq \abs{q_1}-\abs{\ov u-c}-\abs{ \ov v -c}-\abs{u-\ov u} -\abs{v- \ov v}, $$

and using the definitions of $\ov u, \ov v$ and (\ref{generalizedA1eq4}) we get:
\begin{equation}\lab{generalizedA1eq7}
\abs{q_1'} \geq M_0-2\abs z -10\delta-8 r_0-max\{\abs y , \abs{x_2} \}> 3 k \ max\{\abs y, \abs{x_2} \}. 
\end{equation}

We consider  $q_1'=t_1...t_n$, where either  $lab(t_{2i-1})=y^k, \ lab(t_{2i})=x_2^k$ for every $1<2i\leq n-1$ or $lab(t_{2i-1})=x_2^k, \  lab(t_{2i})=y^k$ for every $1<2i\leq n-1$. By the estimate on $q_1'$ above we have that  $n \geq 4$.  

Now we use (\ref{generalizedA1eq2}), (\ref{generalizedA1eq3}) and (\ref{generalizedA1eq7}) to obtain:
$$\abs{q'_2} = \abs{u'-v'} \geq \abs{\ov u'- \ov v'} -\abs{u'- \ov u} - \abs{v'- \ov v} \geq 
\abs{\ov u- \ov v} -\abs{\ov u'- \ov u}-\abs{\ov v'- \ov v}-\abs{u'- \ov u} - \abs{v'- \ov v}$$
 $$\geq M_0-16\delta-2\abs z -10\delta-8 r_0-max\{\abs y , \abs{x_2} \}>3 k \ max\{\abs y, \abs{x_2} \}.$$
We consider  $(q_2')^{-1}=t'_1...t'_{n'}$, where   $lab(t_{2i})=y^{k'}, lab(t_{2i+1})=x_2^{k'}$ for every $1<2i<n'-1$ or $lab(t_{2i})=x_2^{k'}, lab(t_{2i+1})=y^{k'}$ for every $1<2i<\abs{n'}-1$ where $k=\pm k'$. By the estimate on $q_1'$ above, $n' \geq 4$.

 We can now apply lemma \ref{lemma 2.5} to the closed path $p_1'q_1'p_2'q_2'$ with upper bound on $\abs{p_i'}$ provided by (\ref{generalizedA1eq6}) and obtain  a constant $m_0$ such that for every $k \geq m_0$  the paths $t_2$ and $t_3$ are compatible with $t_i'$ and $t_{i+1}'$ respectively (for some unique $i$). Let us denote for convenience $lab(t_2)=W_2 ^k$, $lab(t_3)=W_3 ^k$, $lab(t_i')=\ov W _i ^{k'}$, $lab(t_3)=\ov W _{i+1} ^{k'}$, where the sets $\{W_2,W_3\}$, $\{\ov W_i,\ov W_{i+1}\}$ and $\{x,y\}$ are all equal. Lemma \ref{lemma 2.5} also provides that there exist compatibility paths $v_2$ and $v_3$ with labels $V_2, V_3$ such that:
\begin{align} \lab{38}
	V_2^{-1}W_2^rV_2=\ov W _i^s, \ 
	V_3^{-1}W_3^{r'}V_3=\ov W _{i+1}^{s'},
\end{align}
for some $r,s,r',s'>0$. Because $x_2$ and $y$ are non-commensurable,  the equalities (\ref{38}) are only possible if  $W_2 \equiv \ov W  _i^{\pm 1}$ and $W_3 \equiv \ov W  _{i+1}^{\pm 1}$. Moreover, one of the exponents is positive because $y$ is not conjugate with $y^{-1}$ and thus  $lab (q_2)^{-1}=(y^kx_2^k)^{ m}$ and $W_2 \equiv \ov W  _i$ and $W_3 \equiv \ov W  _{i+1}$. Now by definition of compatible paths we have that $V_2 \in E(W_2)$ and $V_3 \in E(W_3)$. Consider a path $v$ connecting the terminal vertex of $t_2$ with the terminal vertex of $t'_i$. We  also consider a pair of paths $\ov{q_1}v_2 \ov {q_2}$ and $\ov{q_3}v_3 \ov {q_4}$ each of which has the same initial and terminal vertices as the path $v$. Reading off their labels provides the following inequalities in $G$: 
$$lab(v)=W_2 ^{s_1} V_2  W_2 ^{s_2}= W_3 ^{s_3} V_3  W_3 ^{s_4}$$  

for some exponents $s_i  \in \mathbb Z.$  Hence $lab(v) \in E(W_2) \cap E(W_3)=E(x_2) \cap E(y)=E(G)$. We obtained that $z=lab(p_1)$ is equal to either $(y^kx_2^k)^{s'} lab(v)^{-1} (y^kx_2^k)^{-s''}$  or $(y^kx_2^k)^{s'} y^k lab(v)^{-1} ((y^kx_2^k)^{s''}y^k)^{-1}= (y^kx_2^k)^{s'} lab(v)^{-1} (y^kx_2^k)^{-s''}$ for some non-negative numbers $s',s''.$ In both cases $z \in E(G) \times \langle (y^kx_2^k) \rangle.$ 

We obtained that $E(g) \cap \HH=\{e\}$  for $g=y^k x_2^k$,  and now the lemma follows from theorem \ref{quasiconvex and isomorphic}.$\Box$ \medskip

 {\bf Proof of theorem \ref{generalizedA3}} 
 (1) The sufficiency is provided by theorem \ref{quasiconvex and isomorphic}. Assume that there exist an element $x$ of infinite order and $t \neq 0$  such that $\langle  \HH, x^t \rangle \cong  \HH*\langle x^t \rangle$. Take $h \in E(x) \cap \HH$, then there exist $n \neq 0$ and $n'= \pm n$ such that $h^{-1}x^nhx^{n'}=e$ in $G$. Thus $h^{-1}x^{tn}hx^{tn'}=e$ in $G$ which imlpies $h=e$.
 
 (2) The sufficiency follows from lemma \ref{generalizedA1}. To show the necessity it is enough to notice that if the element $x$ satisfying part (1) exists then $E(G) \cap \HH \leq E(x) \cap \HH =\{e \}$. 
 
 (3) Denote the subgroup $E(G) \cap (\HH* \langle x^t \rangle)$ by $K$, it is a finite subgroup in $\HH*\langle x^t \rangle$. By Kurosh subgroup theorem, $K$ is conjugate to a subgroup in $\HH$. On the other hand $K$ is normal in $\HH*\langle x^t \rangle$ and thus $K < \HH$, i.e. $E(G) \cap \HH* \langle x^t \rangle \leq E(G) \cap \HH=\{ e\}$
   
{\bf Proof of corollary \ref{generalizedA2}}

 Consider a canonical homomorphism $\phi: G \rightarrow \ov{G}=G/E(G)$. It is clear that  $E(\ov G)=\{e \}$: the subgroup $E(\ov G)$ is finite normal, hence the subgroup  $\phi^{-1}(E(\ov G))$ is finite normal and thus $\phi^{-1}(E(\ov G)) \leq E(G)$. Homomorphism $\phi$ is a quasi-isometry because $E(G)$ is finite. Thus $\ov \HH =\phi(\HH )$ is quasiconvex of infinite index in $\ov G$. We can apply lemma \ref{generalizedA1}  to $ \ov \HH , \ov G $ and find some $ \ov y$ such that $ E( \ov y )$ is infinite cyclic and obtain the isomorphism $\langle \HH , y \rangle  \cong \ov  \HH *   \langle \ov  y \rangle.$ Consider some preimage $y$ of  $\ov y$. 
$\phi^{-1} (\langle \ov \HH, \ov y \rangle)= \langle \HH \cdot E(G), y \rangle \cong  \HH \cdot E(G)  *_{E(G)}\langle y,E(G) \rangle. $
$\Box$ \medskip
 
 \section {Acknowledgments} The author is grateful to prof. A. Olshanskiy for suggesting the topic of this paper and valuable discussions. I would also like to thank Denis  Osin for his thoughtful suggestions. The author is also greatful to Ashot Minasyan for referring him  to papers \cite{M-P}, \cite{Min} after the first version of this paper was published on Arxiv.org.


\begin{thebibliography}{xx}

\bibitem[Arzh]{Arzh}   \emph{ G. N. Arzhantseva}, On Quasiconvex Subgroups of Word Hyperbolic Groups. Geometriae Dedicata 87: 191--208, 2001.
 
\bibitem[AL]{AL} \emph{G. N. Arzhantseva, I. G. Lysenok}  Growth tightness for word hyperbolic groups. Math. Z. 241 (2002), no. 3, 597–611. 

 
\bibitem[BO]{BO} \emph{Y. Bahturin, A. Olshanskiy}
Actions of maximal growth. Proc. Lond. Math. Soc. (3) 101 (2010), no. 1, 27–72. 

\bibitem[CDP]{CDP} \emph{M. Coornaert, T. Delzant et A. Papadopoulos}, Geometrie et Theorie des Groupes, les Groupes Hyperboliques de Gromov, Lecture Notes in Math., n 1441, Springer, 1990.


\bibitem[Ghys]{Ghys}   \emph{E. Ghys, P. de la Harpe}, Sur Les Groupes Huperboliques D'apres Mikhael Gromov. Birkh\"{a}user, 1990.

\bibitem[GMRS]{GMRS} \emph{R. Gitik, M. Mitra, E. Rips, M. Sageev},
Widths of subgroups. Trans. Amer. Math. Soc. 350 (1998), no. 1, 321–329.

\bibitem[Gre]{Gre}\emph{F.P. Greenleaf}, Invariant Means on Topological Groups and Their Applications, Van Nostrand Reinhold (1969).

\bibitem[Gro]{Gro}  \emph{M. Gromov}, Hyperbolic Groups. Essays in group theory, MSRI Publications,
        Springer, 1987.



\bibitem[IvOl]{IvOl}  \emph {S. V. Ivanov, A. Yu. Olshanskiy }, Hyperbolic groups and their quotients of bounded exponents. Trans. Amer. Math. Soc. 348 (1996), pp 2091-2138.  

\bibitem[Mack]{Mack}   \emph{T. P. Mack}, Quasiconvex Subgroups and Nets in Hyperbolic Groups. PhD Thesis, CalTech 2006. http://thesis.library.caltech.edu/2461/1/thesis.pdf


\bibitem[M-P]{M-P} \emph{Eduardo Martínez-Pedroza}
On residual properties of word hyperbolic groups, GEOMETRIAE DEDICATA, to appear, http://www.springerlink.com/content/y70t6026j8486241.

\bibitem[Min]{Min} \emph{A. Minasyan}  On residual properties of word hyperbolic groups, J. of Group Theory 9 (2006), No. 5, pp. 695-714.




\bibitem[LSch]{LSch}  \emph{M. Lyndon,P. Schupp}, Combinatorial group theory.
        Springer, 1977.



\bibitem[Olsh93]{Olsh93}  \emph {A. Yu. Olshanskiy}, On residualizing homomorphisms and G-subgroups of hyperbolic groups. IJAC, Vol.3, No.4, 1993.






\bibitem[Sta] {Sta} \emph{J. R. Stallings}, Topology of finite graphs. Invent. Math., 71 (1983), 551--565.


\end{thebibliography}
\end{document}